\documentclass[11pt,twoside]{article}
\usepackage{fancyhdr}
\usepackage[colorlinks,citecolor=blue,urlcolor=blue,linkcolor=blue,bookmarks=false]{hyperref}
\usepackage{amsfonts,epsfig,graphicx}
\usepackage{afterpage}
\usepackage{amsmath,amssymb,amsthm} 
\usepackage{fullpage}
\usepackage[T1]{fontenc} 
\usepackage{epsf} 
\usepackage{graphics} 
\usepackage{amsfonts,amsmath}
\usepackage[sort,numbers]{natbib} 
\usepackage{psfrag,xspace}
\usepackage{color,etoolbox}
\usepackage{footmisc}

\theoremstyle{plain}
\newtheorem{theorem}{Theorem}

\newtheorem{lemma}{Lemma}

\newtheorem{fact}{Fact}

\newcommand{\inprod}[2]{\ensuremath{\langle #1 , \, #2 \rangle}}

\newcommand{\plugincov}{\widehat{\psi}_{\text{pi}}^{\text{cov}}}
\newcommand{\firstcov}{\widehat{\psi}_{\text{fo}}^{\text{cov}}}

\newcommand{\truecov}{\psi^{\text{cov}}}
\newcommand{\minimaxcov}{\mathfrak{M}_n^{\text{cov}}(\mathcal{G}(r_n,s_n))}

\newcommand{\plugindens}{\widehat{T}_{\text{pi}}^{\text{f}}}
\newcommand{\firstdens}{\widehat{T}_{\text{fo}}^{\text{f}}}

\newcommand{\truedens}{T^{\text{f}}(f^*)}
\newcommand{\minimaxdens}{\mathfrak{M}_n^{\text{f}}(\mathcal{F}(r_n))}

\newcommand{\pluginseq}{\widehat{Q}_{\text{pi}}^{\theta}}
\newcommand{\firstseq}{\widehat{Q}_{\text{fo}}^{\theta}}
\newcommand{\higherseq}{\widehat{Q}_{\text{ho}}^{\theta}}

\newcommand{\adaseq}{\widehat{Q}_{\text{ad}}^{\theta}}
\newcommand{\trueseq}{Q(\theta^*)}
\newcommand{\minimaxseq}{\mathfrak{M}^{\theta}_n(\Theta(r_n))}

\begin{document}

\begin{center} {\LARGE{\bf{The Fundamental Limits of \\
\vspace{.3cm}
Structure-Agnostic Functional Estimation}}}
\\

\vspace*{.3in}

{\large{
\begin{tabular}{ccccc}
Sivaraman Balakrishnan$^\dagger$, Edward Kennedy$^\dagger$ and Larry Wasserman$^\dagger$ \\
\end{tabular}

\vspace*{.1in}

\begin{tabular}{ccc}
Department of Statistics and Data Science$^{\dagger}$ \\
\end{tabular}

\begin{tabular}{c}
Carnegie Mellon University, \\
Pittsburgh, PA 15213.
\end{tabular}

\vspace*{.2in}

\begin{tabular}{c}
{\texttt{\{siva,edward,larry\}@stat.cmu.edu}}
\end{tabular}
}}

\vspace*{.2in}

\today
\vspace*{.2in}
\begin{abstract}
Many recent developments in causal inference, and functional estimation problems more generally, have been motivated by the 
fact that classical one-step (first-order) debiasing methods, or their more recent sample-split 
double machine-learning avatars, can outperform plugin estimators under surprisingly weak conditions. These first-order corrections
improve on plugin estimators in a black-box fashion, and consequently are often used in conjunction with powerful off-the-shelf estimation methods. 
On the other hand, these first-order methods are 
provably suboptimal in a minimax sense for functional estimation when the nuisance functions live in H\"{o}lder-type function spaces. 
This suboptimality of first-order debiasing has motivated the development of ``higher-order'' debiasing methods~\cite{robins2008higher,van2014higher,birge1995estimation,bickel1988estimating}. 
The resulting estimators 
are, in some cases, provably optimal over H\"{o}lder-type spaces, but in sharp contrast to first-order estimators, both the estimators which 
are minimax-optimal and their analyses are crucially tied to properties of the underlying function space. 
Along a similar vein, some work~\cite{van2014targeted,benkeser2017doubly,dukes2021doubly} has considered $\sqrt{n}$-consistent estimation of causal effects under weaker conditions
than those required by first-order methods, once again relying on higher-order debiasing. 
More recent work in this area has focused on attempting to weaken the dependence of these higher-order 
estimators on the underlying nuisance function spaces, to make the resulting estimators and theory more robust. A central focus
has been to try to make higher-order methods compatible with black-box nuisance estimators. 

In this paper we investigate the fundamental limits of structure-agnostic functional estimation, where relatively weak conditions are 
placed on the underlying nuisance functions. We show that there is a strong sense
in which \emph{existing first-order methods are optimal}. Particularly, we show that for 
several canonical integral functionals of interest it is impossible to improve on first-order 
estimators without making further, strong structural assumptions.
We achieve this goal by providing a formalization of the problem 
of functional estimation with black-box nuisance function estimates, and deriving minimax lower bounds for this problem. Our results highlight some clear tradeoffs in functional estimation -- if we wish to remain agnostic to the underlying nuisance 
function spaces, impose only high-level rate conditions, and maintain compatibility with black-box nuisance estimators then first-order methods 
are optimal. When we have a better understanding of the structure of the underlying nuisance functions then carefully constructed higher-order estimators can outperform 
first-order estimators.
\end{abstract}
\end{center}

\section{Introduction}Statistical modeling often begins by hypothesizing that the data at hand are sampled from a potentially complex, high-dimensional distribution, and the goal in a variety of applications is
not to estimate the distribution itself, but rather to estimate some informative \emph{functional} of the sampling distribution. 
Such functional estimation problems arise naturally in causal inference where under various identification assumptions, causal estimands are expressed as \emph{functionals} of the observed data generating distribution. 
One of the main 
challenges in causal inference is to design statistically efficient functional estimates, while remaining as agnostic as possible to the structure of the sampling distribution (the so-called nuisance component). 
Beyond causal inference, functional estimation problems arise routinely in machine learning \cite{singh2020estimating,samy2014influence}, information theory \cite{paninski2003estimation,jiao2015minimax,wu2016minimax,koza1987statistical}, theoretical computer science \cite{valiant2011estimating} and other fields. 

Recent research in machine learning has led to the development of powerful prediction methods, which perform surprisingly well despite the complexity of the underlying 
prediction tasks as well as the high-dimensionality of the covariates \cite{lecun2015deep}. Consequently, a flurry of research in causal inference \cite{chernozhukov2022automatic,chernozhukov2017orthogonal} has aimed to leverage these 
prediction methods to estimate causal estimands. At the heart of these works 
is the observation that classical one-step/first-order bias-corrected estimators of many important functionals can be constructed to leverage essentially arbitrary initial estimates of the nuisance functions. 
These first-order estimators interact with the nuisance function estimates in a black-box manner, improving on na\"{i}ve plugin estimates by shrinking their bias, but otherwise inheriting their structure-agnostic strengths. Essentially, if 
we are able to construct nuisance function estimates with small error (i.e. typically solve a prediction or density estimation problem well) then the one-step estimator produces an accurate functional estimate. 
An important aspect of this procedure is that we don't need to be able to quantify the precise structure in the nuisance functions that allows us to solve the nuisance function estimation problem well, we simply inherit fast rates
of convergence when we are able to do so. We refer to estimators of this type as \emph{structure agnostic}. Structure agnostic functional estimates are particularly powerful because modern machine learning algorithms 
are in practice able to solve complex prediction tasks with high-dimensional covariates with high accuracy, but we are still far from being able to accurately quantify from a theoretical perspective the precise structures which enable this. 
An important question, one which we aim to formalize and answer in this paper, is: \emph{what are the fundamental limits of structure agnostic functional estimation?}

In many cases, if we can further ensure
that the nuisance estimates converge at a faster than $n^{1/4}$-rate the resulting one-step estimators achieve fast $\sqrt{n}$-rates of convergence, attain semiparametric efficiency bounds, and allow for straightforward
inference~\cite{bickel1993efficient,chernozhukov2017orthogonal,kennedy2022semiparametric}. These ideas are particularly powerful when used together with sample-splitting and cross-fitting, where the nuisance functions 
are estimated on one half of the data, the functional is estimated on the held-out data, and the roles are reversed and the two resulting estimates are averaged to regain efficiency. 

Despite their many strengths, one-step estimators are known to be far from minimax-optimal for many non-parametric functional estimation problems over smoothness classes, even when they are based on a minimax-optimal nuisance function estimate. This basic observation dates back to at least the work of \citet{bickel1988estimating} who constructed minimax-optimal estimates of the integral of the square of a density by further debiasing the one-step estimator. More generally, for estimating
smooth integral functionals of a density~\citet{birge1995estimation} proposed a general higher-order debiasing scheme, and developed complementary lower bounds. Their scheme, in combination with ideas from the papers~\cite{laurent1996efficient,kerk1996estimating,tchetgen2008minimax}, yields minimax-optimal estimates for a broad class of smooth integral functionals. For more complex functionals which arise in causal inference, the construction
of higher-order estimators is more involved, and is the main contribution of a more recent line of work \cite{robins2008higher,van2014higher,robins2017minimax, mukherjee2017semiparametric}, with complementary lower bounds appearing in the work of~\citet{robins2009semiparametric}. In these settings, higher-order estimators improve on the one-step estimate in (very) 
low-regularity settings when $\sqrt{n}$-rates are not achievable, and are also able to achieve $\sqrt{n}$-rates in a wider range of (moderately) low-regularity settings. 
Inspired by this latter observation, some work~\cite{van2014targeted,benkeser2017doubly,dukes2021doubly} has considered $\sqrt{n}$-consistent estimation of causal effects under weaker conditions
than those required by first-order methods, once again relying on higher-order debiasing. It is worth noting that these higher-order estimators improve on first-order estimates, and are minimax-optimal in certain settings, but are decidedly not structure
agnostic in the same way that the plugin and first-order functional estimates are\footnote{We note that we sometimes emphasize certain differences between certain classes of estimators, but the classification of estimators
and the boundaries between these classes can be blurry. Part of the motivation of our work is to ground the discussion of the relative merits of different types of estimators in a rigorous minimax framework.\label{footone}}. 

\vspace{.3cm}

\noindent {\bf Our Contributions: } With this background in place we can now briefly summarize our most significant contributions:
\begin{enumerate}
\item In Section~\ref{sec:minimaxsetup} we describe a formal minimax setup aimed at understanding the 
fundamental limits of black-box functional estimation. This minimax framework allows us to frame the discussion of structure-agnostic 
versus structure-aware estimators, and study their relative merits.
\item In Theorem~\ref{thm:mainlb}, we develop consequences 
for estimating three canonical functionals -- the quadratic functional in the Gaussian sequence model, the quadratic functional in the non-parametric density model, and a
mixed bias causal functional (the expected conditional covariance). 
Building on relatively well-understood techniques, in Theorem~\ref{thm:mainub} we give matching upper bounds. 
Taken together these results highlight the impossibility of improving on first-order estimators without making additional structural
assumptions. 
\item We conclude in Section~\ref{sec:discussion} with some discussion of our results, their implications, and some important avenues for future research.
\end{enumerate}

\subsection{Related Work}
Functional estimation problems have a rich history in many different fields and we refer the reader to the works \cite{tsiatis2006semiparametric,tsybakov2009introduction,van2000asymptotic,bickel1993efficient} for a broader
introduction to the subject. We focus in this section on briefly reviewing some lines of work which provided most of the inspiration for our work, and which study 
functional estimation problems in a minimax framework. In our work we present concrete results for three canonical functional estimation problems:
estimating a non-linear functional in the Gaussian sequence model, estimating a non-linear integral functional of a density, and estimating a causal functional (the expected
conditional covariance). 

Functional estimation in the Gaussian sequence model goes back to the work of \citet{ibragimov1985nonparametric} who initiated the study of 
linear functionals in this model. The work of \citet{donoho1990minimax,fan1991estimation,cai2011testing} have considered estimating non-linear functionals in the Gaussian sequence model,
over Sobolev ellipsoids, Besov bodies, $\ell_p$ balls, and hyperrectangles. More recent work \cite{collier2017minimax,comminges2021adaptive} has focused on estimating linear and non-linear
functionals over sparsity classes. Going beyond the Gaussian sequence 
model, estimating functionals of more general parameter vectors, for instance, functionals of regression coefficients has
also been studied in recent work~\citep{carpentier2022estimation,cai2018semisupervised}.  

The estimation of smooth integral functionals of densities was considered by \citet{bickel1988estimating} who studied estimating the integral of the square of a smooth density.  \citet{bickel1988estimating} 
showed that minimax rates for this functional exhibited an ``elbow effect'' -- when the smoothness of the density $\alpha > d/4$ it is possible to attain parametric $\sqrt{n}$-rates but for less smoothness the best 
achievable rate is non-parametric. Their work inspired further work on estimating other smooth integral functionals of densities and regression functions over non-parametric smoothness classes 
\cite{birge1995estimation,tchetgen2008minimax,kerk1996estimating,laurent1996efficient} culminating in a relatively comprehensive minimax theory for these functionals. The work on estimating 
integral functionals
of a smooth density foreshadowed many developments in causal inference: particularly identifying, the sub-optimality of plugin estimates, the improved but minimax sub-optimal performance of one-step corrected
estimates, and finally minimax-optimal estimates constructed via higher-order corrections. 
There is also a large body of work studying the estimation of two-sample integral density functionals~\citep{berrett2023efficient,singh2018minimax,singh2020estimating,krishnamurthy2014nonparametric,manole2024plugin} which include many natural divergence measures between distributions.
Departing from the minimax framework, there are numerous other frameworks in which one could compare estimators. 
These results often provide a complementary picture. For instance, the work of \citet{cattaneo2022average} studies estimators of the quadratic functional via the inferential lens of bootstrap consistency, highlighting 
other tradeoffs between some of the estimators that we study.

Functionals which arise in causal inference typically exhibit more complex structure, often depending on multiple nuisance functions. The work of \citet{robins1994estimation,robins1995semiparametric}
highlighted the so-called double robustness phenomenon, where the one-step corrected estimates exhibited (second-order) bias which depended on the product of the errors of nuisance estimates. More recent work, 
highlights the benefits of sample-splitting and cross-fitting when using first-order estimates~\cite{chernozhukov2018double}, and attempts to characterize more precisely the set of functionals for which the first-order estimate is doubly robust \cite{rotnitzky2019characterization,chernozhukov2022automatic}.
Moving beyond first-order estimates, the work on higher-order influence functions \cite{robins2008higher,robins2017minimax} and the work on complementary minimax lower bounds~\cite{robins2009semiparametric}
has aimed to more completely develop the minimax theory for various important functionals in causal inference, when the nuisance functions are H\"{o}lder smooth.

Our setting and results also bear similarities to 
the local minimax framework which originates in work of \citet{hajek1972local} and \citet{lecam1972limits}. \emph{Asymptotic} local minimax results are part of the foundation of semi-parametric efficiency theory~\cite{van2000asymptotic,bickel1993efficient}, and are well-studied for functionals in the Gaussian sequence model~\cite{bickel1993efficient,tsybakov2009introduction}, various integral density functionals~\cite{birge1995estimation,laurent1996efficient,berrett2019efficient,hall1986powerful,vanes1992estimating}, and functionals which arise in casual inference (see, for instance, the works \cite{tsiatis2006semiparametric,kennedy2022semiparametric} and references therein).

\subsection{Notation}
We will use the notation $\lesssim, \gtrsim$ to denote inequalities which hold up to a universal positive constant, and $\asymp$ to denote an equality which holds up to a universal positive constant.

\vspace{.3cm}

\noindent {\bf Sobolev Ellipsoids: } Some of our results will consider estimation of functionals in the Gaussian sequence model.
We will discuss, for instance, the case when $\Theta$ is a Sobolev ellipsoid, i.e. for some constants $M_1,M_2 > 0$ our parameter $\theta^*$ is in the set:
\begin{align}
\label{eqn:sobolev}
\Theta^s(M_1,M_2) = \left\{\theta: \sum_{j=1}^\infty \theta_j^2 j^{2s/d} \leq M_1, \sum_{j=1}^\infty \theta_j^2 \leq M_2 \right\}.
\end{align}
The parameter $d > 0$ plays the role of the dimension. When a function $f: \mathbb{R}^d \mapsto \mathbb{R}$ belongs to a Sobolev space and has $s$ weak derivatives with finite $\ell_2$ norm, then under certain conditions, the function can be represented in an orthonormal basis with coefficients which belong to the Sobolev ellipsoid $\Theta^s$. 

\vspace{.3cm}

\noindent {\bf H\"{o}lder Functions: } For a function $f: \mathbb{R}^d \mapsto \mathbb{R}$, and a vector $\alpha \in \mathbb{N}^d$ we define,
\begin{align*}
D^{\alpha} f = \frac{\partial^{\|\alpha\|_1} f }{ \partial x_1^{\alpha_1} \ldots \partial x_{d}^{\alpha_d}}. 
\end{align*}
Then letting $\ell = \lceil s \rceil - 1$ we define the H\"{o}lder function class: 
\begin{align*}
\mathcal{H}^s(L) &= \Big\{f: f~\text{is}~\ell~\text{times differentiable}, \\
&~~~~~~~~~ |D^{\alpha} f(x) - D^{\alpha} f(y)| \leq L\|x - y\|_2^{s - \ell}, \\
&~~~~~~~~~\forall (x,y) \in \mathbb{R}^d, \|\alpha\|_1 = \ell, \alpha \in \mathbb{N}^d \Big\}.
\end{align*}
We sometimes refer to both H\"{o}lder functions and Sobolev ellipsoids with the terminology H\"{o}lder-type spaces. 
\section{Background}
We begin by introducing some functional estimation problems. 
We then briefly introduce classical one-step estimators for this problem emphasizing their structure-agnostic nature, before discussing minimax-optimal higher-order estimators.
\subsection{Functional Estimation Problems}

Although our results have 
broader implications, we focus throughout on three important functional estimation problems, for which minimax rates are relatively well-understood. In each case, we briefly summarize
some well-known results. We revisit some generalizations of our results in Section~\ref{sec:discussion}, and highlight some potential avenues for further investigation.

\vspace{.3cm}

\noindent {\bf Quadratic Functionals in the Gaussian Sequence Model: } The main ideas of our work are most clearly understood in the following classical infinite Gaussian 
sequence model. We observe,
\begin{align*}
y_j = \theta_j^* + \epsilon_j,
\end{align*}
where $j \in \{1,2,\ldots\}$, each $\epsilon_j$ is drawn independently with distribution $N(0,1/n)$ and our goal is to estimate the quadratic functional:
\begin{align}
\label{eqn:seq}
Q(\theta^*) = \sum_{j=1}^{\infty} \theta_j^{*2}. 
\end{align}
This functional is a canonical example of a smooth functional and minimax rates for estimation (over Sobolev ellipsoids) go back to a series
of works \cite{bickel1988estimating,laurent1996efficient,donoho1990minimax,fan1991estimation}.

\vspace{.3cm}

\noindent {\bf Smooth Integral Functionals: } 
In this setting, we observe $X_1,\ldots,X_n \sim f^*$, and our goal is to estimate an integral functional,
\begin{align}
\label{eqn:dens}
T(f^*) = \int (f^*(x))^2 dx. 
\end{align}
This can be generalized to the estimation of $T_\varphi(f^*) = \int \varphi(f^*(x)) dx,$ for some $\varphi$ which has continuous second derivative. Under suitable conditions, this general setup 
includes the estimation of familiar information-theoretic quantities (like the entropy), and familiar quantities which arise in non-parametric estimation (like $\ell_p^p$ norms). 
The estimation of the quadratic functional was studied by~\citet{bickel1988estimating} and the more general problem of estimating smooth integral functionals was considered by~\citet{birge1995estimation}.

\vspace{.3cm}

\noindent {\bf Causal Functionals: } To illustrate our ideas we focus on the expected conditional covariance. This functional 
arises in biostatistics and epidemiology in the context of the estimation of the causal effect of a binary treatment \cite{robins2009semiparametric}, and is an important functional in assessing conditional independence. 
Concretely,  we observe samples of the form $\{(Y_1, A_1, X_1), \ldots, (Y_n, A_n,X_n)\}$ drawn i.i.d from a distribution $\mathbb{P}$, where $X_i \in \mathbb{R}^d, A_i \in \{0,1\}, Y_i \in \{0,1\}$. We refer to
$X$ as the covariates, $Y$ as the outcome, and $A$ as the treatment. 
We denote the covariate density $p_X$, and define the regression function and propensity score:
\begin{align*}
\mu^*(x) &= \mathbb{E}[Y | X = x], \\
\pi^*(x) &= \mathbb{E}[A | X = x]. 
\end{align*}
Our goal is to estimate the functional:
\begin{align}
\label{eqn:truecov}
\truecov = \mathbb{E}[\text{cov}(A,Y|X)] = \mathbb{E} [AY] - \int \pi^*(x) \mu^*(x) p_X(x).
\end{align}
The first term is easy to estimate at fast $\sqrt{n}$-rates, and the focus is often on estimating the second term. Even in our 
binary setup the joint distribution over triples $(X,A,Y)$ is not fully specified by the nuisance functions $(\mu^*, \pi^*,p_X)$. The joint distribution 
in the binary setup is fully parametrized by \emph{quadruples} $(\mu^*, \pi^*, \eta^*, p_X)$ where we additionally define,
\begin{align*}
\eta^*(x) = \mathbb{E}(Y | X, A = 1) - \mathbb{E}(Y | X, A = 0).
\end{align*}
Somewhat surprisingly minimax rates for estimating the expected conditional covariance are not understood in full generality when the nuisance functions are H\"{o}lder smooth. 
In the more restricted setting when the covariate density $p_X$ is either known or can be estimated at a sufficiently fast rate, sharp minimax rates
are better understood~\cite{robins2008higher}.

\subsection{The Methodological Approach to Functional Estimation}
The functional estimation problems introduced above are examples of semi-parametric inference problems, with non-parametric nuisance components. 
The first attempt to solve these problems was based on the so-called plugin principle. Using the data we estimate the non-parametric components and then plug them in to obtain estimates
of the functional. For instance, in the problem of estimating the expected conditional covariance, we regress the outcome $Y$ on the covariates $X$ to obtain an estimate $\widehat{\mu}$, 
and regress the treatment $A$ on the covariates $X$ to obtain an estimate $\widehat{\pi}$ and construct the plugin estimate:
\begin{align}
\label{eqn:plugincov}
\plugincov = \frac{1}{n} \sum_{i=1}^n A_i Y_i - \frac{1}{n} \sum_{i=1}^n \widehat{\pi}(X_i) \widehat{\mu}(X_i). 
\end{align}
Similar plugin estimates can be constructed for functionals in the Gaussian sequence model, and for density integral functionals. It is natural to expect that if our estimates $\widehat{\mu}$ and 
$\widehat{\pi}$ are accurate, then the resulting plugin estimate will also be accurate.  
It is worth emphasizing the \emph{structure-agnostic nature} of the plugin estimates, and particularly its compatibility with black-box nuisance function estimates. One can use a powerful machine learning algorithm
(say a random forest, or a deep neural network) to construct estimates of the propensity score and the regression function, and use these to construct accurate functional estimates.

One drawback of plugin estimates is that they inherit bias, and rates of convergence, directly from their nuisance function estimates. This in turn can complicate inference, and has led to the development 
of powerful one-step correction methods which improve on plugin estimates by reducing their bias, improving their rate of convergence, and allowing for valid $\sqrt{n}$-rate inference, even when the nuisance function
estimates converge at a slower than $\sqrt{n}$-rate. Developing and analyzing one-step corrected estimators in a variety of different contexts is the focus of a large body of past work (see the recent review~\cite{kennedy2022semiparametric} and references therein). 
These one-step corrections are at the heart of semi-parametric theory, and are particularly powerful when used in conjunction with sample-splitting (and/or cross-fitting). To illustrate the main idea, suppose 
we consider the expected conditional covariance, and consider the following first-order estimator:
\begin{align}
\label{eqn:firstcov}
\firstcov = \frac{1}{n} \sum_{i=1}^n (A_i - \widehat{\pi}(X_i))(Y_i - \widehat{\mu}(X_i)),
\end{align}
where now we suppose that $\widehat{\pi}$ and $\widehat{\mu}$ are constructed on a separate sample. It is common to construct two estimates (reversing the roles of the two samples) and average them.
This estimator can be viewed as arising from correcting the plugin estimator in~\eqref{eqn:plugincov} by adding to it an estimate of the \emph{influence function} of the target functional.
To build some intuition for the estimate $\firstcov$ one can observe that treating $\widehat{\mu}, \widehat{\pi}$ as fixed (or estimated on a separate sample) we can compute 
the bias (or conditional bias) of $\firstcov$ and observe that it is of second-order,~i.e. 
\begin{align*}
\left| \mathbb{E}[\firstcov] - \truecov \right| &= \Big| \int (\pi^*(x) - \widehat{\pi}(x))  \\ 
&~~~~~~~~(\mu^*(x) - \widehat{\mu}(x)) p_X(x) dx \Big|.
\end{align*}
The estimator $\firstcov$ exhibits the so called double robustness property. Roughly, one can upper bound the error of the estimate $\firstcov$ by a product of errors of the underlying nuisance estimates.
Once again, it is worth emphasizing the structure-agnostic nature of the first-order estimate. As with the plugin estimate, the first-order estimate is agnostic to the nature of the underlying nuisance function estimates. The  
guarantees for the first-order estimator rely on the accuracy of the pilot estimates, but neither the estimator nor its guarantee are tailored to the structure which enabled accurate pilot estimation. This enables
us to use this functional estimate along with black-box machine learning algorithms, which perform well in practice, but do so by exploiting 
structural properties of the underlying nuisance functions that can be difficult to describe mathematically. 

\subsection{Smoothness Classes and The Structural Approach to Functional Estimation}
\label{sec:struct}
Functional estimation problems are also studied from a minimax perspective in order to understand fundamental limits and to construct optimal estimators. Absent any structural 
assumptions consistent functional estimation is impossible and it is classical to impose some structure on the non-parametric components in the form of smoothness assumptions.

In the Gaussian sequence model this amounts to constraining the parameter space by hypothesizing that $\theta^* \in \Theta$. Minimax rates for quadratic functional estimation are well-understood 
for a large variety of constraint sets $\Theta$. In our discussion, we will primarily focus on the case when $\theta^*$ is in the Sobolev ellipsoid $\Theta^s(M_1,M_2)$ in~\eqref{eqn:sobolev}.
In the case of integral functionals of a density, it is common to hypothesize that the nuisance function (the sampling density $f$) belongs to a H\"{o}lder space, i.e. that $f \in \mathcal{H}^s(L)$, and 
that $f$ specifies a valid density $f \geq 0$, 
$\int f(x) dx = 1.$
Finally, for the expected conditional covariance a typical assumption is that $\pi \in \mathcal{H}^\alpha(L_1)$, $0 \leq \pi(x) \leq 1$, and $\mu \in \mathcal{H}^{\beta}(L_2)$. 

Given these structural assumptions, it is then natural to wonder: given a rate-optimal (say in the $\ell_2$-sense) nuisance function estimate, are the resulting plugin or first-order
estimators minimax optimal? The answer to this question is often ``no'' 
and this in turn motivates higher-order estimators.

First-order estimators can be viewed as a linear bias correction of a plugin estimate. These first-order estimators have quadratic bias, i.e. in typical cases the bias decays quadratically in the error of the plugin estimator (see, for instance~\cite{robins2008higher} for a detailed discussion).
Higher-order estimators are constructed, roughly, by subtracting
an estimate of this quadratic bias from the first-order estimate. Higher-order U-statistics are commonly used to estimate parameters which can written as the expected values of functions of pairs (or higher-order tuples) of independent samples. In functional estimation problem the estimate of the bias of the first-order estimator takes the form of a higher-order U-statistic \citep{bickel1988estimating,robins2008higher,birge1995estimation}.
In the Gaussian sequence model, one second-order estimator for $Q(\theta^*)$ is a classical truncated series estimator.
For a truncation threshold $T > 0$ we construct the estimator:
\begin{align}
\label{eqn:hoseq}
\higherseq = \sum_{j=1}^T  y^1_j y^2_j +  \sum_{j = T+1}^{\infty} \left[(y_j^1 + y_j^2) \widehat{\theta}_j - \widehat{\theta}_{j}^2\right], 
\end{align}
where $\widehat{\theta}_j$ denotes the $j$-th coordinate of the pilot estimate $\widehat{\theta}$. We assume for simplicity that we obtain two observations $y^1, y^2$ in the Gaussian sequence model. When this is not the case, one can use the sample-splitting device described in~\cite{robins2006adaptive}, and define:
\begin{align*}
y_j^1 &= y_j + \Phi^{-1}(U_j)/\sqrt{n} \\
y_j^2 &= y_j - \Phi^{-1}(U_j)/\sqrt{n},
\end{align*}
where $U_j$ are independently drawn uniform random variables. The resulting $y_j^1, y_j^2$ are now independent, with means $\theta_j^*$ and variance $2/n$ (i.e. their variances are inflated 
by a factor of $2$).
When $\widehat{\theta}$ is 0, the estimate~\eqref{eqn:hoseq} is an unbiased estimate of $\sum_{j=1}^T \theta_j^{*2}.$ 
In this case, it is well-known~\cite{laurent1996efficient,donoho1990minimax,fan1991estimation}, that when the truncation parameter $T$ is chosen to scale as $n^{2d/(4s + d)}$, $\higherseq$ is a minimax optimal 
estimate of $Q(\theta^*)$ over $\Theta^s$ and is semi-parametrically efficient when $s > d/4$. 

It is important to note that despite the fact the estimator was constructed as a second-order correction to the plugin estimate, its minimaxity over the Sobolev ellipsoid is no longer strongly dependent 
on the choice of the pilot estimate $\widehat{\theta}$, which could simply be taken to be 0.  Rather, the optimality 
of this estimator is closely related to properties, such as decay rate of the coefficients $\theta_j^*$, of the underlying Sobolev ellipsoid. 
In a similar vein, higher-order estimators have been constructed for integral functionals of a density in~\cite{kerk1996estimating,birge1995estimation,tchetgen2008minimax}, and for causal functionals like the expected conditional covariance in~\cite{robins2008higher}. These estimators
exhibit similar properties to the estimator~\eqref{eqn:hoseq}. Demonstrating their minimax optimality (or advantage over plug-in or first-order estimates) often requires carefully leveraging the properties of the underlying nuisance function space.

\vspace{.3cm}

\noindent {\bf Comparing Higher-order Estimators and One-Step Corrections:  } To summarize the discussion so far, it is worth once again contrasting first-order and higher-order estimators to identify
some common themes which hold across the canonical examples and more broadly. First-order estimators are black-box corrections to plugin estimates. Under very weak conditions they improve on plugin 
estimators, and their accuracy depends only on the (squared) errors of the nuisance function estimates. This in turn enhances their compatibility with black-box (machine learning) methods for nuisance function
estimation.  They are often not minimax-optimal over H\"{o}lder-type function spaces, even when used in conjunction with minimax-optimal nuisance function estimates. In contrast, minimax-optimal higher-order
estimators are more carefully tailored to the underlying function space. They typically have a weaker connection to the plugin estimates on which they are based. 
They are analyzed via a more careful understanding of the bias-variance tradeoff in the underlying nuisance function space. They can in many cases yield minimax-optimal functional estimates over
H\"{o}lder-type spaces, even when used with trivial (zero) pilot nuisance function estimates. 

\section{Main Results}
We begin with a description of the black-box minimax setup we will focus on. We then turn our attention to minimax lower bounds in Section~\ref{sec:lbmain}, and briefly provide complementary upper bounds
in Section~\ref{sec:ubmain}.

\subsection{Minimax Functional Estimation in the Black-Box Model}
\label{sec:minimaxsetup}
To fix ideas we first consider estimation of the quadratic functional in the Gaussian sequence model~\eqref{eqn:seq}. As we noted previously, with no structural assumptions placed on $\theta^*$ 
consistent functional estimation is impossible. Rather than impose smoothness assumptions, we model the black-box setting where we construct a pilot estimate on a separate sample. 

More formally, our goal is to estimate $Q(\theta^*)$ and our assumption
on $\theta^*$ is that the pilot estimate $\widehat{\theta}$ is accurate in an $\ell_2$ sense, i.e. that $\theta^* \in \Theta(r_n)$, where:
\begin{align}
\label{eqn:bbcon}
\Theta(r_n) := \left\{ \theta: \|\theta - \widehat{\theta}\|_2^2 \leq r_n \right\},
\end{align}
where the accuracy of the pilot estimate $r_n$ is unknown to the statistician. 

It is important to note that the assumption that the pilot estimate $\widehat{\theta}$ is $r_n$-accurate, imposes an (implicit) structural assumption on $\theta^*$. The strength of this structural 
assumption depends on the (unknown) rate of convergence $r_n$.  \emph{It is precisely
this structural condition that plugin and first-order estimates are tailored to leverage.} It is also worth contrasting this structural assumption with the smoothness assumptions in~\eqref{eqn:sobolev} -- here the structural assumption 
hypothesizes that our favorite nuisance function estimator returns an accurate pilot estimate, but does not further constrain $\theta^*$ to have a particular structure.
Since $r_n$ is unknown to the statistician, estimators constructed in this model are implicitly \emph{adaptive}, i.e. this setting bears similarities to the classical adaptive non-parametric estimation setting where
functions are hypothesized to be smooth but the smoothness parameter is unknown to the statistician.

Absent any additional smoothness assumptions, our goal is to construct a minimax rate-optimal estimate, i.e. an estimate $\widehat{Q}$ such that:
\begin{align}
\label{eqn:mmseq}
&\sup_{\theta^* \in \Theta(r_n)} \mathbb{E} (\widehat{Q} - Q(\theta^*))^2 \asymp \inf_{\widetilde{Q}} \sup_{\theta^* \in \Theta(r_n)} \mathbb{E} (\widetilde{Q} - Q(\theta^*))^2 \\
&~~~~~~~~~~~~~~~~~~~~~~~~=: \minimaxseq, \nonumber
\end{align}
and to study the minimax risk $\minimaxseq$.  

In a similar vein, one can consider minimax estimation of the density functional $T(f^*)$ in~\eqref{eqn:dens}. We assume for simplicity that the densities under consideration 
are uniformly upper bounded by some (large) constant $M > 0$. 
Our goal is to estimate $T(f^*)$ under the constraint that $f^* \in \mathcal{F}(r_n)$:
\begin{align*}
\mathcal{F}(r_n) &:= \Big\{f:  \int (f(x) - \widehat{f}(x))^2 dx \leq r_n, f \geq 0, \\
&~~~~~~~~~~\int f(x) dx = 1,  \|\widehat{f}\|_{\infty}, \|f\|_{\infty} \leq M\Big\}. 
\end{align*}
In this case we define the associated minimax risk as:
\begin{align}
\label{eqn:mmdens}
\minimaxdens := \inf_{\widehat{T}} \sup_{f^* \in \mathcal{F}(r_n)} \mathbb{E} (\widehat{T} - T(f^*))^2.
\end{align}
For the expected conditional covariance we are given two pilot estimates $\widehat{\mu}$ and $\widehat{\pi}$. In order to construct higher-order estimators we
might further assume that we are given a third pilot estimate $\widehat{p}_X$ of the
covariate density. To simplify our presentation of minimax lower bounds, 
we consider the case when the 
covariate density $p_X$ is uniform on $[0,1]^d$.
We also assume that $\widehat{\mu}$ and $\widehat{\pi}$ are bounded
away from 0 and 1 on $[0,1]^d$. This latter restriction can be eliminated via a perturbation argument similar to the one used in~Supplementary Section~B.2 for the integral of
the squared density. Although these restrictions ease the construction of minimax lower bounds, the upper bounds in Theorem~\ref{thm:mainub} hold without these 
restrictions. It is important to note that designing accurate pilot estimates without assuming a lower bound on the density is significantly more challenging. Recent work~\cite{berrett2019efficient,berrett2023efficient} has explored this issue, characterizing how the tails of the sampling distribution influences the complexity of functional estimation problems.

With this setup in place our goal is to estimate $\truecov$ in~\eqref{eqn:truecov}, under the following constraints on $(\mu^*, \pi^*,\eta^*,p_X^*) \in \mathcal{G}(r_n, s_n)$:
\begin{align*}
\mathcal{G}&(r_n, s_n) := \Big\{(\mu,\pi, \eta, p_X): \text{supp}(X) = [0,1]^d, \\
&p_X = \text{unif}[0,1]^d, \int (\mu(x) - \widehat{\mu}(x))^2 p_X(x) dx \leq r_n, \\
& \int (\pi(x) - \widehat{\pi}(x))^2 p_X(x) dx \leq s_n, 0 \leq \pi(x),\mu(x) \leq 1, \\
&1-\varepsilon \geq \widehat{\pi}(x),  \widehat{\mu}(x) \geq \varepsilon > 0, \text{for}~x \in [0,1]^d \Big\}.
\end{align*}
Here we tacitly suppress the dependence of $\mathcal{G}$ on $\varepsilon$ which we 
treat as a universal constant. 
We define the associated minimax risk as:
\begin{align}
\label{eqn:mmcov}
\minimaxcov := \inf_{\widehat{\psi}} \sup_{(\mu^*,\pi^*) \in \mathcal{G}(r_n,s_n)} \mathbb{E} (\widehat{\psi} - \truecov)^2.
\end{align}
With this setup in place, our goal is to understand the fundamental limits on structure-agnostic  
functional estimation by providing upper and lower bounds on the minimax risks in \eqref{eqn:mmseq},~\eqref{eqn:mmdens} and~\eqref{eqn:mmcov}. 

\subsubsection{Interpreting the Minimax Setup}
\label{sec:interpretations}
The minimax problems described in this section require some care to interpret.
We describe briefly some interpretations focusing again on the Gaussian sequence model:

\vspace{.3cm}

\noindent {\bf Sample-Splitting and the Conditional Viewpoint: } When analyzing black-box sample-splitting based functional estimators, it is natural to take a conditional viewpoint to remain judicious
in the modeling assumptions we impose. In this viewpoint, we hypothesize
that for some function $r_{n,\delta}$ our nuisance estimates are $r_{n,\delta}$-accurate with probability at least $1-\delta$, and then proceed to analyze a functional estimate constructed on a separate sample. This perspective
is for instance explicitly adopted in the work~\cite{foster2019orthogonal}, and is implicit in a long series of past work~\cite{bickel1993efficient,bickel1988estimating,chernozhukov2022automatic,chernozhukov2018double}.

To complement this conditional viewpoint with minimax lower bounds, one can aim to understand the fundamental limitations of the second stage of this two-stage estimator construction (the first-stage corresponds to the 
well-studied problem of function estimation). This problem is at the heart of our proposed minimax setup.

In contrast to the traditional setting, where $\theta^*$ is fixed, and $\widehat{\theta}$ is random, in our lower bounds we treat both as fixed. 
To model more closely the sample-splitting based functional estimation paradigm, we might split the data into two sets $\mathcal{D}_0$ and $\mathcal{D}_1$, and construct a pilot estimate $\widehat{\theta}$ on $\mathcal{D}_0$.
We would then relax the constraint set in~\eqref{eqn:bbcon} to be the random set:
\begin{align*}
\Theta(r_n,\delta) = \begin{cases}
\left\{ \theta: \|\theta - \widehat{\theta}\|_2^2 \leq r_n \right\}~\text{with probability}~1 - \delta \\
~~~~~~~~~~~~~~~~~~~~~~~~~~~~~~~~~~~\text{independent of $\mathcal{D}_1$}, \\
\left\{ \theta: \sum_{j=1}^\infty \theta_j^2 \leq M_2 \right\}~\text{otherwise,}
\end{cases}
\end{align*}
where with probability $\delta$ the nuisance parameter $\theta^*$ is only constrained to have bounded norm, and the pilot estimate is uninformative. 
Lower bounds over this random constraint set, can directly be obtained from lower bounds over the set~\eqref{eqn:bbcon} at the cost of some additional notational burden. This is because 
when the nuisance parameter is only constrained to have bounded norm, non-trivial functional estimation is impossible, and on this event the best estimate of the functional is the trivial estimate of $0$.

\vspace{.3cm}

\noindent  {\bf Unstructured Local Minimax Lower Bounds: } An alternative way to interpret our problem setting, is as a type of local minimax setup. The pilot estimate $\widehat{\theta}$, and the local radius $r_n$,
define a local estimation problem via the constraints in~\eqref{eqn:bbcon}. The minimax rates we study are thus quantifying the difficulty of functional estimation, locally around the pilot estimate, with an (otherwise) 
unconstrained nuisance parameter.

Despite the similarity in spirit, the setup and goals are quite different from classical local minimax problems. 
We don't make strong assumptions on the sequence $r_n$, the only constraints in our problem are locality constraints as opposed to locality and smoothness constraints, and we adopt a non-asymptotic view. The most significant difference is that 
we localize the parameter space around the pilot estimate $\widehat{\theta}$ and not the true parameter $\theta^*$ since our goal is not to 
elicit local minimax rates in the neighborhood of the true parameter $\theta^*$, but rather to understand the fundamental limits of black-box functional estimation.

\subsection{Lower Bounds}
\label{sec:lbmain}
With the minimax setup introduced in the previous section we are now equipped to state our lower bounds on the minimax risks in~\eqref{eqn:mmseq},~\eqref{eqn:mmdens} and~\eqref{eqn:mmcov}.
\begin{theorem}
\label{thm:mainlb}
The minimax risks in~\eqref{eqn:mmseq},~\eqref{eqn:mmdens} and~\eqref{eqn:mmcov} are lower bounded as:
\begin{enumerate}
\item {\bf Quadratic Functional in the Gaussian Sequence Model: }
\begin{align*}
\minimaxseq \gtrsim r_n^2 + \|\widehat{\theta}\|^2_2 \min\left\{r_n , \frac{1}{n} \right\}. 
\end{align*}
\item {\bf Quadratic Density Integral Functional: }
\begin{align*}
&\minimaxdens \gtrsim r_n^2 + \\
&\left[ \int \widehat{f}(x)^3 dx - \Big(\int \widehat{f}(x)^2\Big)^2 dx \right] \min\left\{r_n, \frac{1}{n} \right\}.
\end{align*}
\item {\bf Expected Conditional Covariance:} 
\begin{align*}
\minimaxcov &\gtrsim   r_n \times s_n + \frac{1}{n}. 
\end{align*}
\end{enumerate}
\end{theorem}
\noindent We prove this result in the Supplementary Section~B. Our proofs are based on a well-understood recipe. We reduce the problem of lower bounding the minimax risk 
of functional 
estimation to lower bounding the risk (the sum of Type I and II errors) in an appropriate hypothesis testing problem (roughly of distinguishing if the functional is large or small). If the 
null and alternate 
are difficult to distinguish, we obtain a lower bound on the minimax risk. To lower bound the minimax hypothesis 
testing error we carefully construct
priors on the composite null and composite alternate and show that the resulting mixture distributions are difficult to distinguish by lower 
bounding the error
of the (optimal) likelihood ratio test.
At a more technical level, we use classical ideas from \citet{ingster2003nonparametric} for the Gaussian sequence model, from~\citet{balakrishnan2019hypothesis}
for the density integral functional, and from \citet{robins2009semiparametric} for the expected conditional covariance, in order to upper bound an appropriate
divergence measure between the two mixture distributions.

It is interesting to first focus on the case when $r_n,s_n \gg 1/n$, since this is the typical case. In this case, the second term in each of the lower bounds corresponds (in rate)
to semi-parametric efficiency lower bounds. These lower bounds are determined by the \emph{variance} of the estimated influence function, and decay to 0 at the standard parametric rate.
A large fraction of semi-parametric theory focuses on imposing conditions on the nuisance functions under which the \emph{bias} term (of order $r_n^2$ or $r_n \times s_n$) above
decays to zero faster than the parametric rate, at which point semi-parametric efficient estimation and inference are possible. On the other hand, the more recent literature on higher-order
estimators~\cite{robins2008higher,van2014higher,birge1995estimation,bickel1988estimating,van2014targeted,benkeser2017doubly,dukes2021doubly}, has aimed at reducing the bias term to either obtain $\sqrt{n}$-rates under weaker assumptions on the nuisance functions, or in order to obtain (slower than $\sqrt{n}$) 
minimax-optimal rates when the nuisance functions are H\"{o}lder smooth. 
The main import of Theorem~\ref{thm:mainlb} is that under only the assumption that the pilot estimates are accurate in an $\ell_2$ sense, and in 
the absence of further smoothness assumptions, \emph{no further bias reduction is possible.} As we explore further in Theorem~\ref{thm:mainub}, first-order estimates
are minimax-optimal, and achieve the limits of structure-agnostic functional estimation.

When the condition that $r_n,s_n \gg 1/n$ is violated 
the first-stage pilot estimates are super-accurate, i.e. are more accurate than the (fixed-dimensional) parametric rate, and the plugin functional estimate
is already minimax optimal. We address this situation in more detail in the Supplementary Section~D, where we construct a functional estimate 
which adapts between the plugin and first-order estimates paying only a small statistical price, and prove a matching lower bound which highlights
the fundamental limits of adaptivity in this setup. When the pilot estimate is super-accurate, it can be beneficial to use the plugin estimator and ignore the data entirely. More broadly, in some functional estimation problems, it can be the case that seemingly natural `oracle' estimators, for instance estimators which use extra knowledge of nuisance functions in a natural way, can be improved by estimators which do not use this extra knowledge \cite{su2023when,hirano2003efficient,berrett2023efficient,robins1992estimating}.

\subsection{Upper Bounds}
\label{sec:ubmain}
In this section, we develop upper bounds on the minimax risk by analyzing plugin and first-order estimators.
Our analyses of these estimators are elementary, and plugin and first-order estimators have been analyzed much more generally in past work, although often from an asymptotic perspective~\cite{van2000asymptotic}.
Our simple non-asymptotic analysis enables a more direct comparison with lower bounds from Theorem~\ref{thm:mainlb}. 

To set the stage we first formally define the plugin and first-order estimates. 
The plugin and first-order estimates for the quadratic functionals in the Gaussian sequence model and for the density integral functional are:
\begin{alignat*}{2}
\pluginseq &= \|\widehat{\theta}\|_2^2 & \plugindens &= \int (\widehat{f}(x))^2 dx \\
\firstseq &=  2 \inprod{y}{\widehat{\theta}} - \|\widehat{\theta}\|_2^2~~~~~~& \firstdens &= \frac{2}{n} \sum_{i=1}^n \widehat{f}(X_i) - \int (\widehat{f}(x))^2 dx. 
\end{alignat*}
We also recall the definitions of the plugin and first-order estimates of the expected conditional covariance in~\eqref{eqn:plugincov} and~\eqref{eqn:firstcov}, 
and the higher-order estimate in the Gaussian sequence model in~\eqref{eqn:hoseq}.
Having introduced the plugin and first-order estimates of our three canonical functionals we have the following theorem:
\begin{theorem}
\label{thm:mainub}
The estimators described above have the following guarantees:
\begin{enumerate}
\item {\bf Quadratic Functional in the Gaussian Sequence Model: } 
\begin{align*}
|\pluginseq - \trueseq|^2 &\lesssim r_n^2 + r_n \|\widehat{\theta}\|_2^2 \\
\mathbb{E} (\firstseq - \trueseq)^2 &\lesssim r_n^2 + \frac{\|\widehat{\theta}\|_2^2}{n} \\
\mathbb{E} (\higherseq - \trueseq)^2 &\lesssim \left[\sum_{j = T+1}^{\infty} (\widehat{\theta}_j - \theta_j^*)^2 \right]^2 + \frac{\|\widehat{\theta}\|_2^2}{n} + \frac{T}{n^2} \\
&\lesssim r_n^2 + \frac{r_n}{n} + \frac{\|\widehat{\theta}\|_2^2}{n} + \frac{T}{n^2}. 
\end{align*}
\item {\bf Quadratic Density Integral Functional: }
\begin{align*}
|\plugindens - \truedens|^2 &\lesssim r_n^2 + r_n \int \widehat{f}(x)^2 dx \\
\mathbb{E} (\firstdens- \truedens)^2 &\lesssim r_n^2 + \frac{\text{var}(\widehat{f}(X))}{n}.
\end{align*}
\item {\bf Expected Conditional Covariance:} 
\begin{align*}
\mathbb{E} (\plugincov - \truecov)^2 &\lesssim r_n \times s_n + r_n + s_n + \frac{1}{n} \\
\mathbb{E} (\firstcov- \truecov)^2 &\lesssim  r_n \times s_n + \frac{1}{n}. 
\end{align*}
\end{enumerate}
\end{theorem}
\noindent Once again we initially focus our discussion on the case when $r_n, s_n \gg \frac{1}{n}$, which is the typical setting, in which case the first-order estimates outperform the plugin estimates.
We design
an adaptive estimate for the quadratic functional in the sequence model, which selects between the plugin and first-order estimate in a data-driven manner, and achieves close to the oracle risk, in Supplementary Section~D.
For each of our functionals, focusing on terms which only depend on $n$ and $r_n, s_n$, the first-order estimate achieves a maximum risk which matches the minimax lower bounds of Theorem~\ref{thm:mainlb}, i.e. the first-order estimates
are minimax optimal in our setting. 

The higher-order estimator in the Gaussian sequence model depends on a truncation parameter $T$. 
The higher-order estimator can have smaller bias than the first-order estimator since it unbiasedly estimates $\sum_{j=1}^T \theta_j^{*2}$, and combines this with a first-order estimator of $\sum_{j=T+1}^\infty \theta_j^{*2}$.
When $T$ is set small relative to $n$, then this estimator incurs only a modest amount of additional variance. To achieve minimax-optimality over a Sobolev ellipsoid in the low-regularity regime when $s < d/4$, the truncation
parameter $T$ needs to be chosen larger than $n$, in a careful way, to balance the reduction in bias with the inflation in variance relative to the first-order estimator. In our structure-agnostic model, it is impossible to 
guarantee that the higher-order estimator has meaningfully lower bias than the first-order estimator, and consequently it is impossible to guarantee that the higher-order estimator
improves on the first-order estimator. 
Higher-order estimators for the quadratic density functional and for the expected conditional covariance are more involved to describe~\cite{birge1995estimation,robins2008higher} but share the same
qualitative features. 
This in turn highlights that in our formalization of the black-box, structure-agnostic functional estimation problem, where we are unwilling to assume
more than access to a potentially accurate black-box prediction algorithm, it is impossible to improve on first-order estimators in a minimax sense. 
Higher-order estimates can only improve on first-order estimates when additional structural assumptions are imposed and exploited.

We note in passing that there are some slight differences between the constant factors in the upper and lower bounds for the quadratic density functional. This is due to the mismatch between the $\ell_2^2$ 
distance and the squared Hellinger distance. Intuitively these terms of the lower bound are determined by the modulus of continuity of the quadratic functional over a squared Hellinger neighborhood
around $\widehat{f}$~\cite{donoho91geometrizing}. On the other hand, the upper bounds are determined by the modulus of continuity over an $\ell_2^2$ neighborhood. These coincide when we assume that the densities under
consideration are both upper and lower bounded by universal constants, but can otherwise differ. 

Finally, for the expected conditional covariance, we observe that 
the first-order estimator is minimax-optimal for any choice of $r_n, s_n$. This is because even when our pilot estimates are super accurate, i.e. for instance the true $\pi^*$ and $\mu^*$ are known,
we still need to estimate the term $\mathbb{E}[AY]$ to construct our functional estimate. This in turn leads to an unavoidable $\mathcal{O}(1/n)$ term in the MSE. In contrast, for the quadratic density and sequence functionals 
the plugin estimator is a deterministic function of the nuisance estimates and incurs no additional variance, and can dominate the first-order estimator when the pilot estimates are super-accurate. 

\section{Discussion and Extensions}
\label{sec:discussion}
In this work, we introduced a minimax framework for reasoning about two-stage structure-agnostic functional estimation methods. We developed consequences
for estimating three canonical functionals -- the quadratic functional in the Gaussian sequence model, the quadratic functional in the non-parametric density model, and a
mixed bias causal functional (the expected conditional covariance). 

By focusing on concrete examples, we have given results for particular estimators which are canonical plugin and first-order estimators, but have avoided giving precise general definitions for these classes of 
estimators. This is by design, as we noted in Footnote~\ref{footone} the distinctions between these classes of estimators can be blurry. For instance, the work of \citet{newey2018cross} and \citet{gine2008simple}
show that certain carefully undersmoothed plugin-type estimators can perform similarly to higher-order estimators, and inherit both their strengths and weaknesses.

There are several possible extensions of our results. From a technical standpoint, for the quadratic density integral functional and the expected conditional covariance, our upper and lower bounds do not match in their dependence
on the constants defining the model classes. Tightening these discrepancies could be interesting and might also suggest alternative formulations of the structure-agnostic minimax rate which build on pilot estimators which are accurate in other metrics than the ones we have used in this paper.
Minimax theory is well-understood for more general smooth functionals in both the density model and in the Gaussian sequence model, 
and this theory mirrors closely results for the quadratic functionals in these models. We expect our main results will continue to hold for more general smooth functionals. For causal functionals, beyond the expected 
conditional covariance, minimax lower bounds are known only in a few problems \cite{robins2009semiparametric} and we expect our results will extend to cover these functionals as well. A more ambitious
extension would aim to cover larger classes of functionals for which first-order estimators are well-understood to have desirable properties~\cite{rotnitzky2019characterization,chernozhukov2022automatic}. However,
the likelihood structure in each of these statistical models is quite different, which poses some challenges to developing a unified theory of lower bounds. Minimax rates, over appropriate smoothness classes, have also been studied
for certain local integral functionals which arise in causal inference and in non-parametric regression~\cite{shen2020optimal,wang2008effect,cai2009variance,kennedy2022minimax}, and it would be also be interesting to understand the fundamental limits of 
structure-agnostic estimation in these problems.

Finally, it is worth emphasizing that our lower bounds do not preclude estimators which improve in some restricted 
ways on first-order estimators. For instance, it is possible in some cases to construct adaptive 
estimators which perform nearly as well
as the first-order estimator in the absence of any additional structure, but improve on the first-order estimator
when the nuisance functions have additional smoothness structure. Such estimates are developed, for instance, in a testing context in the work of \citet{liu2020nearly}. 
On the other hand, our results do show that if one aims to improve on first-order estimators in a general minimax sense, this improvement is only possible at the expense of adding
further assumptions, i.e. there are limits to what can be achieved by higher-order estimators without additional structural assumptions. 
Developing a comprehensive understanding of 
adaptive estimators, which adapt between smoothness classes and structure-agnosticity, and understanding their fundamental limitations, could be an interesting avenue for future research. 

\section*{Acknowledgements}
The authors are grateful to Jamie Robins for several helpful discussions, and for many inspiring conversations. The authors are also grateful to the anonymous reviewers for their insightful comments and careful reading, which significantly strengthened the paper.

\bibliographystyle{abbrvnat}
\bibliography{bibliography}

\appendix

\section{Lower Bound Preliminaries}
In this section, we collect some well-known technical facts that will aid the proof of Theorem~\ref{thm:mainlb}. 
Suppose that our goal is to estimate a functional $T(P)$, given samples from $P \in \mathcal{P}$. 
We first recall a standard construction (see for instance Theorem 2.14 in~\cite{tsybakov2009introduction}) for obtaining
lower bounds in functional estimation problems.

We construct two prior distributions, 
$\pi_0$ and $\pi_1$ on $\mathcal{P}$, which induce two distributions $Q_0$ and $Q_1$ where for any measurable set $A$: 
\begin{alignat*}{2}
Q_0(A) &= \int P^n(A) d\pi_0(P),~~~&\text{and}~~~ Q_1(A) &= \int P^n(A) d\pi_1(P). 
\end{alignat*}
We further ensure that our functional takes sufficiently different values under each of the prior distributions, i.e.:
\begin{align*}
\pi_0(\{P: T(P) \leq c\}) = 1,~~~&~~~\pi_1(\{P: T(P) \geq c + 2s\}) = 1.
\end{align*}
We denote by $\chi^2(P,Q)$ the chi-squared divergence between $P$ and $Q$, and by $H^2(P,Q)$ the squared Hellinger distance between $P$ and $Q$ (see Section 2.4 of the book \cite{tsybakov2009introduction} for formal definitions).
Then we have the following result:
\begin{lemma}
\label{lem:fuzzy}
Suppose that $\chi^2(Q_0, Q_1) \leq \alpha < \infty$, then 
\begin{align*}
\inf_{\widehat{T}} \sup_{P \in \mathcal{P}} \mathbb{E}(\widehat{T} - T(P))^2 \geq s^2 \max \Big( \frac{1}{4} \exp(-\alpha), \frac{1 - \sqrt{\alpha/2}}{2} \Big). 
\end{align*}
If $H^2(Q_0, Q_1) \leq \alpha < 2$, then:
\begin{align*}
\inf_{\widehat{T}} \sup_{P \in \mathcal{P}} \mathbb{E}(\widehat{T} - T(P))^2 \geq s^2 \frac{1 - \sqrt{\alpha(1 - \alpha/4)}}{2}. 
\end{align*}
\end{lemma}
\noindent The proof of this lemma follows immediately from Theorem 2.15 in~\cite{tsybakov2009introduction}, and an application of Markov's inequality to obtain in-expectation bounds. 
The main takeaway is simply that if we can construct $Q_0, Q_1$ as above, ensuring that the functional is separated by at least $2s$, and ensuring that $Q_0$ and $Q_1$ have 
$\chi^2$ divergence or Hellinger distance upper
bounded by a sufficiently small constant, then we obtain lower bounds on the minimax risk of order $s^2$. 

\vspace{.3cm}

\noindent We frequently use the following well-known fact. 
\begin{fact}
\label{fac:one}
Given two measures $p$, $q$, the squared Hellinger distance between their $n$-fold products can be upper bounded as $H^2(p^n, q^n) \leq n H^2(p,q)$.
\end{fact}

Suppose we consider the Gaussian sequence model and let $Q_0 = N(\theta, 1/n)$. 
We then define $Q_1$ in the following way. 
We select any $d$ coordinates, let $\mathcal{I}$ denote the selected indices. Fix an $\epsilon > 0$. Let $N := 2^d$ 
and let $\{u_1,\ldots,u_{N}\}$ denote the collection of all vectors with entries $\{+\epsilon, -\epsilon\}$. 
For any $v \in \mathbb{R}^d$ we 
denote by
$\theta_{v}$ the vector which perturbs $\theta$ by adding $v$ to the indices in $\mathcal{I}$. 
Define $Q_1^\epsilon$ to be the following mixture:
\begin{align*}
Q_1^\epsilon = \frac{1}{N} \sum_{i=1}^{N} N(\theta_{u_i},1/n).
\end{align*}
That is, the mixture distribution $Q_1^\epsilon$ is obtained by perturbing the coordinates in $\mathcal{I}$ by $\pm \epsilon$ uniformly at random.
The following result is well-known:
\begin{lemma}
\label{lem:chisq}
For any $\theta$, and $\mathcal{I}$, suppose $\epsilon \leq 1/n$, then:
\begin{align*}
\chi^2(Q_0, Q_1^\epsilon) \leq \exp(d n^2 \epsilon^4) - 1.
\end{align*}
Furthermore, if $d n^2 \epsilon^4 \leq 1$, then, 
\begin{align*}
\chi^2(Q_0, Q_1^\epsilon) \leq 2 d n^2 \epsilon^4.
\end{align*}
\end{lemma}
\noindent We include a short proof for completeness in~Appendix~\ref{app:chisq}.
We also note the following simple bound on the $\chi^2$ distance between two Gaussians.
\begin{lemma}
\label{lem:simple}
Suppose $P = N(\theta, I/n)$, $Q = N(\widetilde{\theta},I/n)$, then 
\begin{align*}
\chi^2(P,Q) = \exp(  n \|\theta - \widetilde{\theta}\|_2^2) - 1. 
\end{align*}
Furthermore, if $n \|\theta - \widetilde{\theta}\|_2^2 \leq 1$, then,
\begin{align*}
\chi^2(P,Q) \leq 2 n \|\theta - \widetilde{\theta}\|_2^2.
\end{align*}
\end{lemma}

We will use the following constrained risk inequality as a consequence of Theorem 1 of \citet{brown1996constrained}. Though we only need a result for the squared loss, the 
form we use is most directly deduced from Corollary~2 of \cite{duchi2018constrained}. This inequality is a useful tool for proving lower bounds for adaptive estimators.

\begin{lemma}
\label{lem:cri}
Fix any $\theta$, and define $D_1 := N(\theta, I/n)$ and $D_2 := N((1 + \alpha/\|\theta\|_2)\theta, I/n)$. 
Let $P_1$ denote either $D_1$ or $D_2$, and let $P_2$ denote the other distribution. 
Suppose we have an estimator $\widehat{Q}$ for which,
\begin{align*}
\mathbb{E}_{P_1}(\widehat{Q} - Q(P_1))^2 \leq \beta^2, 
\end{align*}
then 
\begin{align*}
\mathbb{E}_{P_2}(\widehat{Q} - Q(P_2))^2 \geq \left[\alpha^2 + 2 \alpha \|\theta\|_2 - \beta \exp(n \alpha^2/2) \right]^2_+ .
\end{align*}
\end{lemma}

\subsection{Proof of Lemma~\ref{lem:chisq}}
\label{app:chisq}
The proof follows from a direct computation. We note that,
\begin{align*}
\chi^2(Q_0, Q_1^\epsilon) &= \mathbb{E}_{Q_0} \left( \frac{Q_1^\epsilon}{Q_0}\right)^2 - 1 \\
&= \mathbb{E}_{Q_0} \prod_{i=1}^d \left[ \frac{ \frac{1}{2} \left[\exp\left( - n \frac{(y_i - \theta_i - \epsilon)^2}{2} \right) + \exp\left( - n \frac{(y_i - \theta_i + \epsilon)^2}{2} \right)\right]}{\exp(- n(y_i - \theta_i)^2/2)} \right]^2- 1 \\
&\stackrel{\text{(i)}}{=} \mathbb{E}_{Q_0} \prod_{i=1}^d \frac{1}{4} \exp(-n \epsilon^2) \left[ \exp( 2 \sqrt{n} \epsilon Z_i) + \exp( - 2 \sqrt{n} \epsilon Z_i) + 1\right] - 1 \\
&\stackrel{\text{(ii)}}{=} \prod_{i=1}^d \frac{1}{2} \left[ \exp(n \epsilon^2) + \exp(-n \epsilon^2) \right] - 1 \\
&\stackrel{\text{(iii)}}{\leq} \prod_{i=1}^d \exp(n^2 \epsilon^4) - 1 \\
&= \exp(dn^2 \epsilon^4) - 1. 
\end{align*}
In (i), $Z_i$ denotes a standard Gaussian random variable, (ii) uses the fact that $\mathbb{E}(\exp(tZ_i)) = \exp(t^2/2)$, 
and (iii) uses the fact $\text{cosh}(x) \leq 1 + x^2 \leq \exp(x^2)$ for $x \in [0,1]$.

\section{Proof of Theorem~\ref{thm:mainlb} }
\label{app:lower}
We prove each of the three lower bounds in turn.
\subsection{Lower Bounds for Quadratic Functional in the Gaussian Sequence Model} 
\label{app:seqlb}
Our goal is to prove minimax lower bounds for estimating the quadratic functional in  the GSM. In particular, 
we'd like to understand lower bounds on:
\begin{align*}
\inf_{\widehat{Q}} \sup_{\theta: \|\theta - \widehat{\theta}\|_2^2 \leq r_n} \mathbb{E}(\widehat{Q} - Q(\theta))^2. 
\end{align*}
Our lower bounds will be a consequence of Lemma~\ref{lem:fuzzy} with various choices of the priors. 

\noindent {\bf LB 1: } $r_n \geq \frac{1}{n}, r_n^2 \lesssim \frac{\|\widehat{\theta}\|_2^2}{n}$. In this case, we construct $Q_0 = N(\widehat{\theta}, I/n)$ and $Q_1 = N((1- \nu) \widehat{\theta}, I/n)$ where $0 \leq \nu \leq 1$.
Then by Lemma~\ref{lem:simple}, we have that $\chi^2(Q_1,Q_0) \leq \alpha$ if $2 n \nu^2 \|\widehat{\theta}\|_2^2 \leq \alpha$, so we select 
\begin{align*}
\nu^2 = \frac{\alpha}{2n \|\widehat{\theta}\|_2^2},
\end{align*}
to ensure this. Since $r_n^2 \lesssim \|\widehat{\theta}\|_2^2/n$, and $r_n \geq 1/n$, we observe that $\nu < 1$. The distance, $\|\widehat{\theta} - (1 - \nu)\widehat{\theta}\|_2^2 = \nu^2 \|\widehat{\theta}\|_2^2 = \alpha/2n \leq 1/n$ so the two priors are supported on the set $\|\widehat{\theta} - \theta\|_2^2 \leq r_n$ 
as desired.

The functional separation under the two priors is:
\begin{align*}
2s &:= \|\widehat{\theta}\|_2^2 - (1 - \nu)^2 \|\widehat{\theta}\|_2^2 \\
&= (2\nu - \nu^2) \|\widehat{\theta}\|_2^2 \\
&\geq \nu  \|\widehat{\theta}\|_2^2.
\end{align*}
So applying Lemma~\ref{lem:fuzzy}, the minimax rate is at least $\nu^2 \|\widehat{\theta}\|_2^4$, i.e. $\frac{\|\widehat{\theta}\|^2_2}{n}$ as desired. 

\noindent {\bf LB 2: } $r_n \leq \frac{1}{n}$. In this case, we use the pair of distributions, $Q_0 = N(\widehat{\theta}, I/n)$ and $Q_1 = N(\widehat{\theta} + \sqrt{r_n} \widehat{\theta}/\|\widehat{\theta}\|_2, I/n)$.
By Lemma~\ref{lem:simple}, we have that $\chi^2(Q_1,Q_0) \leq \alpha$ if $2n r_n \leq \alpha$. 

The functional separation in this case:
\begin{align*}
2s &:= \|\widehat{\theta}\|_2^2 \left(1 + \frac{\sqrt{r_n}}{\|\widehat{\theta}\|_2} \right)^2 - \|\widehat{\theta}\|_2^2 \\
&= r_n + 2\sqrt{r_n}\|\widehat{\theta}\|_2\\
&\geq 2 \sqrt{r_n}\|\widehat{\theta}\|_2.
\end{align*}
So the minimax rate is at least $r_n \|\widehat{\theta}\|_2^2$ as desired.

\noindent {\bf LB 3: } Finally we show that $r_n^2$ is a lower bound. To see this we observe that without loss of generality we can assume that 
$\widehat{\theta}$ has finite norm since if it did not have finite norm we have already shown that the minimax
rate is infinite. As a consequence, for any $\epsilon > 0$, and for any finite integer $d > 0$, we can find $d$ indices, denoted as $\mathcal{I}$, such that $|\widehat{\theta}_j| \leq \epsilon/4$, for $j \in \mathcal{I}$. 

We construct two mixtures $Q_0$ and $Q_1^\epsilon$ in the following way: we set $Q_0 = N(\widehat{\theta}, I/n)$ and $Q_1^{\epsilon}$ as described in the setup to Lemma~\ref{lem:chisq}, by perturbing each coordinate
in $\mathcal{I}$ by $\pm \epsilon$ uniformly at random. 
Lemma~\ref{lem:chisq} shows that the $\chi^2$ distance is at most $d n^2 \epsilon^4$. On the other hand the functional separation is,
\begin{align*}
2s &:= \left[\|\widehat{\theta}\|_2^2 - \sum_{j \in \mathcal{I}} \widehat{\theta}_j^2 + \sum_{j \in \mathcal{I}} (\widehat{\theta}_j \pm \epsilon)^2 \right] - \|\widehat{\theta}\|_2^2  \\
&= \sum_{j \in \mathcal{I}} \epsilon^2 \pm 2 \widehat{\theta}_j \epsilon \\
&\geq \frac{1}{2} \sum_{j \in \mathcal{I}} \epsilon^2 = \frac{d\epsilon^2}{2}. 
\end{align*}
On the other hand, the distribution $Q_1^\epsilon$ is supported on parameters $\widetilde{\theta}$ such that, $\|\widehat{\theta} - \widetilde{\theta}\|_2^2 = d\epsilon^2$. It remains to prescribe choices for $d, \epsilon$ 
to ensure that our constraints are satisfied. We choose $\epsilon = \sqrt{\alpha} \min\{1/(nr_n), 1/n\}$, and $d = r_n/\epsilon^2$. This ensures that the $\chi^2$ distance is at most $\alpha$, and as a consequence of Lemma~\ref{lem:fuzzy}
we obtain a minimax lower bound of order $r_n^2$ as desired.

\subsection{Lower Bounds for the Integral of the Squared Density} 
\label{app:lowerdens}
We prove lower bounds corresponding to each term in our bound separately. 

{\bf LB 1: } $r_n \gtrsim 1/n$.
We use Le Cam's two-point method. Define, $p_1 := \widehat{f}$ and for a sufficiently small $\varepsilon > 0$ define,
\begin{align*}
p_2(x) := p_1(x) \left( 1 + \varepsilon p_1(x) - \varepsilon \theta \right), 
\end{align*}
where $\theta := \int \widehat{f}(x)^2 dx$. Since, the densities under consideration are bounded by $M > 0$, we observe that for sufficiently small $\varepsilon > 0$, $p_2$ is a valid
density. Now, we observe that the Hellinger distance between these densities,
\begin{align*}
H^2(p_1, p_2) &\leq \int \frac{(p_1(x) - p_2(x))^2}{p_1(x)} dx = \int \varepsilon^2 p_1(x) (p_1(x) - \theta)^2 dx = \varepsilon^2 \left( \int p_1(x)^3 dx - \theta^2\right) := \varepsilon^2 \gamma.
\end{align*}
On the other hand, the functional separation is,
\begin{align*}
\int p_2^2(x) dx - \int p_1^2(x) dx &= 2\varepsilon \left( \int p_1(x)^3 dx - \theta^2\right) + \varepsilon^2 \int (p_1^2(x) - \theta p_1(x))^2 dx \\
&\geq 2\varepsilon \left( \int p_1(x)^3 dx - \theta^2\right) = 2 \varepsilon \gamma.
\end{align*}
Finally, we note that the $\ell_2^2$ distance is upper bounded as,
\begin{align*}
\int (p_1(x) - p_2(x))^2 dx =  \int \varepsilon^2 p_1^2(x) (p_1(x) - \theta)^2 dx \lesssim \varepsilon^2 \gamma, 
\end{align*}
using the fact that $\|p_1\|_{\infty} \leq M$. Now, we set, $\varepsilon$ such that, $\varepsilon^2 \gamma = \alpha/n$, for a sufficiently small constant $\alpha > 0$. We note that this construction is valid since $r_n \gtrsim 1/n$. Then the squared Hellinger 
distance between the $n$-fold product measures is at most $\alpha$ using Fact~\ref{fac:one}, and Lemma~\ref{lem:fuzzy} yields a lower bound of order $\gamma/n$ as claimed.

{\bf LB 2: } $r_n \lesssim 1/n$.
In this case, we follow the same construction as above except we choose $\varepsilon$ such that, $\varepsilon^2 \gamma = \alpha r_n$, for a sufficiently small constant $\alpha > 0$. Then 
Lemma~\ref{lem:fuzzy} yields a lower bound of order $\gamma r_n$ as claimed.

{\bf LB 3: } Lower Bound of Order $r_n^2$.
To set the stage we derive a result on the chi-squared distance between a density, and a perturbed counterpart. Our result and proof are
largely inspired by the proof of Lemma 4.4 in~\citet{balakrishnan2019hypothesis}, which in turn generalizes a result of~\citet{ingster1997adaptive}.

Our construction involves a set $S \subseteq \mathbb{R}^d$, which, for a given $m$, we divide into $2m$ disjoint sets of equal volume, and pair these sets together into $m$ (disjoint) pairs $\{(A_1,B_1), \ldots, (A_m,B_m)\}$. Now, suppose that we construct:
\begin{align*}
p_{\lambda}(x) = p_0(x) + \frac{h}{\sqrt{\text{vol}(A_1)}} \sum_{j=1}^m \left[ \lambda_j \left[\mathbb{I}(x \in A_j) - \mathbb{I}(x \in B_j)\right]\right],
\end{align*}
where $0 \leq h/\sqrt{\text{vol}(A_1)} \leq \inf_{x \in S} p_0(x)$, will be chosen appropriately in the sequel, $\lambda_j \in \{-1,+1\}$ will be chosen uniformly at random, and $\text{vol}(A)$ denotes the Lebesgue measure
of the set $A$. We note that by our choice of $h$, $p_{\lambda}$ is a valid density, i.e. 
is positive everywhere and integrates to 1.
We then have the following result:
\begin{lemma}
\label{lem:denschisq}
Define $Q$ to be the mixture obtained from choosing $\lambda$ uniformly at random from $\{-1,+1\}^m$. If we satisfy the conditions:
\begin{enumerate}
\item $h/\sqrt{\text{vol}(A_1)} \leq \inf_{x \in S} p_0(x)$
\item $nh^2/(\inf_{x \in S} p_0(x)) \leq 1$
\item $\frac{mn^2h^4}{(\inf_{x \in S} p_0(x))^2} \leq \ln(1 + \alpha),$
\end{enumerate}
then we have that,
\begin{align*}
\chi^2(P_0^n, Q^n) \leq \alpha.
\end{align*}
\end{lemma}
\noindent Taking this result as given we can now complete the proof of our result. By our assumption, $\|\widehat{f}\|_{\infty} \leq M$, so we have that $\int \widehat{f}^2 < M$. Given an arbitrary pilot density $\widehat{f}$, and any $\epsilon > 0$ 
we can find a set $S$ of volume $1$ such that, 
$\widehat{f}(x) \leq \epsilon/4$ on $S$. We choose $\epsilon = \min\{\sqrt{r_n}/4\sqrt{1 + M},1\}.$

Now, we construct an intermediate density $\widetilde{p} = (1-\epsilon)\widehat{f} + \epsilon p_U$, where $p_U$ is the density of the uniform distribution on $S$. 
Now, we note that,
\begin{align*}
\int (\widetilde{p}(x) - \widehat{f}(x))^2 dx = \epsilon^2 \int (\widehat{f} - p_U)^2 dx \leq \epsilon^2 + \epsilon^2 \int \widehat{f}(x)^2 dx,
\end{align*}
so for our choice of $\epsilon$ 
we obtain that,
$\int (\widetilde{p}(x) - \widehat{f}(x))^2 dx \leq r_n/4.$ We also note that $\|\widetilde{p}\|_{\infty} \leq M$, since $\|\widetilde{p}\|_{\infty} \leq \max\{\|\widehat{f}\|_{\infty}, 5\epsilon/4\}.$ Now, on the set $S$, $\widetilde{p}$ is lower bounded as $\inf_x \widetilde{p}(x) \geq \epsilon$. 

We construct perturbations $p_\lambda$ as described above centered around the density $\widetilde{p}$, with $h = \epsilon/\sqrt{2m}$ (for an $m$ we will choose in the sequel). 
We use Lemma~\ref{lem:denschisq}
with $h = \epsilon/\sqrt{2m}$, and choose $m$ large enough to ensure that, $n \epsilon/m \leq 1$, and $\frac{n^2\epsilon^2}{m} \leq \ln(1 + \alpha)$. Each of the densities $p_{\lambda}$ are uniformly upper bounded by
a similar argument to the one for $\widetilde{p}$. Then as a consequence, we obtain that 
$\chi^2(\widetilde{P}^n, Q^n) \leq \alpha$, where $Q$ is the distribution defined in Lemma~\ref{lem:denschisq}.
We note that,
\begin{align*}
\int (p_{\lambda} - \widehat{f})^2 dx \leq 2 \int (p_{\lambda} - \widetilde{p})^2 dx +  2 \int (\widetilde{p} - \widehat{f})^2 dx \leq 2 \epsilon^2 + r_n/2 \leq r_n.
\end{align*}

\noindent Finally, the functional separation is:
\begin{align*}
\int p_{\lambda}^2 - \int \widetilde{p}^2  = \int_S p_{\lambda}^2 - \int \widetilde{p}^2 \geq \frac{1}{2} (2\epsilon)^2 - \frac{25}{16} \epsilon^2 \geq \frac{7}{16} \epsilon^2 \gtrsim \min\{r_n/(1 + M),1\},
\end{align*}
as desired.

\subsubsection{Proof of Lemma~\ref{lem:denschisq}}
We first bound the expected squared likelihood ratio, which in turn implies a bound on the chi-squared distance. Let us denote $\widetilde{h} := h/\sqrt{\text{vol}(A_1)}$.
Suppose that we observe $\{Z_1,\ldots,Z_n\}$, then the likelihood ratio:
\begin{align*}
W_n(Z_1,\ldots,Z_n) = \frac{1}{2^m} \sum_{\lambda \in \{-1,+1\}^m} \prod_{i=1}^n \frac{p_{\lambda}(Z_i)}{p_0(Z_i)}, 
\end{align*}
and 
\begin{align*}
W_n^2(Z_1,\ldots,Z_n) &= \frac{1}{2^{2m}} \sum_{\lambda \in \{-1,+1\}^m} \sum_{\nu \in \{-1,+1\}^m} \prod_{i=1}^n \frac{p_{\lambda}(Z_i) p_{\nu}(Z_i)}{p_0^2(Z_i)} \\
&=  \frac{1}{2^{2m}} \sum_{\lambda \in \{-1,+1\}^m} \sum_{\nu \in \{-1,+1\}^m} \prod_{i=1}^n \left(1 + \frac{\widetilde{h} \sum_{j=1}^m \left[ \lambda_j \left[\mathbb{I}(Z_i \in A_j) - \mathbb{I}(Z_i \in B_j)\right]\right]}{p_0(Z_i)}\right) \times \\
&~~~~\left(1 + \frac{\widetilde{h} \sum_{j=1}^m \left[ \nu_j \left[\mathbb{I}(Z_i \in A_j) - \mathbb{I}(Z_i \in B_j)\right]\right]}{p_0(Z_i)}\right).
\end{align*}
Using the fact that the supports of the sets $\{(A_1,B_1),\ldots,(A_m,B_m)\}$ are all disjoint, taking the expected value over $Z_1,\ldots,Z_n$ and using their independence we obtain,
\begin{align*}
\mathbb{E} [W_n^2(Z_1,\ldots,Z_n)] &=  \frac{1}{2^{2m}} \sum_{\lambda \in \{-1,+1\}^m} \sum_{\nu \in \{-1,+1\}^m} \left( 1 + \widetilde{h}^2 \sum_{j=1}^m \lambda_j \nu_j a_j \right)^n,
\end{align*}
where 
\begin{align*}
a_j = \int_{A_j \cup B_j} \frac{1}{p_0(x)} dx. 
\end{align*}
From this we can write, 
\begin{align*}
\mathbb{E} [W_n^2(Z_1,\ldots,Z_n)] &\leq \mathbb{E}_{\lambda, \nu} \exp( n \widetilde{h}^2 \sum_{j=1}^m \lambda_j \nu_j a_j ) = \prod_{j=1}^m \text{cosh}(n \widetilde{h}^2 a_j).
\end{align*}
Now, using the fact that $\text{cosh}(x) \leq 1 + x^2 \leq \exp(x^2)$ for $x \leq 1$, we obtain that,
\begin{align*}
\mathbb{E} [W_n^2(Z_1,\ldots,Z_n)] &\leq \exp(n^2 \widetilde{h}^4 \sum_{j=1}^m a_j^2) \leq \exp\left(\frac{m n^2 h^4}{(\inf_{x \in S} p_0(x))^2} \right),
\end{align*}
provided that, $n\widetilde{h}^2 a_j \leq 1$ for $j \in \{1,\ldots,m\}$. 
Finally, we note that the $\chi^2$ distance is simply, $\mathbb{E} [W_n^2(Z_1,\ldots,Z_n)] - 1$, and so,
\begin{align*}
\chi^2(P_0^n, Q^n) \leq \alpha, 
\end{align*}
when $mn^2h^4/(\inf_{x \in S} p_0(x))^2 \leq \ln(1 + \alpha)$, as claimed.

\subsection{Lower Bounds for the Expected Conditional Covariance}
Now, we turn our attention to the expected conditional covariance. Our proof in this case is inspired by that of Theorem 4.1 of \citet{robins2009semiparametric}.
Following their work, we construct our lower bound in the case when $Y \in \{0,1\}$. 
In this binary setup the joint distribution over triples $(X,A,Y)$ can be parametrized by quadruples $(\mu, \pi, \eta, p_X)$.
Recall our statistical model, 
\begin{align*}
\mathcal{G}&(r_n, s_n) := \Big\{(\mu,\pi, \eta, p_X): \text{supp}(X) = [0,1]^d, p_X = \text{unif}[0,1]^d, \int (\mu(x) - \widehat{\mu}(x))^2 p_X(x) dx \leq r_n, \\
& \int (\pi(x) - \widehat{\pi}(x))^2 p_X(x) dx \leq s_n, 0 \leq \pi(x),\mu(x) \leq 1, 1-\varepsilon \geq \widehat{\pi}(x),  \widehat{\mu}(x) \geq \varepsilon > 0, \text{for}~x \in [0,1]^d \Big\}.
\end{align*}
In our proof, we will often suppress dependence on the (universal) constant $\varepsilon > 0$. 

\vspace{.3cm}

\noindent {\bf LB 1: } To begin with we show a lower bound of order $1/n$, even when the nuisance functions
$\pi, \mu$ are known exactly, i.e. when $r_n = s_n = 0$. 
The likelihood in our model is given as:
\begin{align*}
p(X,A,Y) &= p_X(X)\pi(X)^A(1 - \pi(X))^{1-A} (\mu(X) + (1-\pi(X))\eta(X))^{AY} (\mu(X) - \pi(X) \eta(X))^{(1-A)Y} \times \\
&~~~(1 - \mu(X) - (1 - \pi(X))\eta(X))^{A(1-Y)} (1 - \mu(X) + \pi(X) \eta(X))^{(1-A)(1-Y)}. 
\end{align*}
We define a pair of distributions $p_1, p_2$, defined by quadruples $(\widehat{\mu}, \widehat{\pi}, 0, 1)$ and $(\widehat{\mu}, \widehat{\pi}, \zeta, 1)$, for some sufficiently small $\zeta > 0$. 
In particular, we will choose $\zeta \leq \varepsilon/2$. 
It is straighforward to verify that the functional separation between these distributions is, 
\begin{align*}
\truecov(p_1) - \truecov(p_2) = \zeta \mathbb{E}[\widehat{\pi}^2(X)] \gtrsim \zeta. 
\end{align*}
On the other hand, the Hellinger distance between $p_1$ and $p_2$ can be calculated by noting:
\begin{align*}
p_1(X,A,Y) &= \widehat{\pi}(X)^A (1 - \widehat{\pi}(X))^{1-A} \widehat{\mu}(X)^Y (1 - \widehat{\mu}(X))^{1-Y}, \\
p_2(X,A,Y) &= \widehat{\pi}(X)^A(1 - \widehat{\pi}(X))^{1-A} (\widehat{\mu}(X) + (1-\widehat{\pi}(X)) \zeta)^{AY} (\widehat{\mu}(X) - \widehat{\pi}(X) \zeta)^{(1-A)Y} \times \\
&~~~(1 - \widehat{\mu}(X) - (1 - \widehat{\pi}(X)) \zeta)^{A(1-Y)} (1 - \widehat{\mu}(X) + \widehat{\pi}(X) \zeta)^{(1-A)(1-Y)}. 
\end{align*}
For any $(X,A,Y)$ we note that, $p_1(X, A, Y) \gtrsim \varepsilon^2$, and $|p_1(X,A,Y) - p_2(X,A,Y)| \lesssim \zeta$ from which we obtain that, 
the squared Hellinger distance, 
\begin{align*}
H^2(p_1, p_2) \lesssim \zeta^2. 
\end{align*}
Choosing $\zeta \lesssim 1/\sqrt{n}$, using Fact~\ref{fac:one} and applying Lemma~\ref{lem:fuzzy} then yields a lower bound of order $1/n$ as desired.

\vspace{.3cm}

The remainder of the analysis is devoted to proving a lower bound of order $r_n \times s_n$. 
We remark that in our setup this lower bound is simpler to derive that lower bounds for the expected conditional covariance under smoothness
conditions \cite{robins2009semiparametric}, since we are not constrained by potentially different smoothnesses
of the propensity score and outcome regression function. 
At a more technical level, it is not necessary in our setup to construct mixtures under both 
the null and alternate, and the bound on the Hellinger distance we need is essentially due to \citet{birge1995estimation}.
It is also not necessary to use different constructions depending on which of the propensity score or outcome regression
is more difficult to estimate. 

\vspace{.3cm}

\noindent {\bf LB 2: } In this case, we are aiming for a lower bound of $r_n \times s_n$, in the setting when $r_n \times s_n \gtrsim 1/n$. 
Consequently, we focus on lower bounds for the estimation of 
\begin{align*}
\psi = \int \pi(x) \mu(x) p_X(x) dx,
\end{align*}
since as we noted earlier the remaining term in the expected conditional covariance can be estimated at fast $\sqrt{n}$-rates. 

\noindent Fix an integer $2m$, and denote by ${B}_1,\ldots, {B}_{2m}$ be $2m$ 
translates of the cube $(2m)^{-1/d} [0,1/2]^d$ which are disjoint, and contained in $[0,1]^d$, and let the bottom left corners of these cubes be $x_1,\ldots,x_{2m}$. 

We now define, 
\begin{align*}
\pi_\lambda(x) &= \widehat{\pi}(x) + \frac{h_1}{\sqrt{\text{vol}(B_1)}\widehat{\mu}(x)} \sum_{j=1}^{m} \lambda_j [\mathbb{I}(x \in B_{2j}(x)) - \mathbb{I}(x \in B_{2j-1}(x))], \\
\mu_\lambda(x) &= \widehat{\mu}(x) + \frac{h_2}{\sqrt{\text{vol}(B_1)} \widehat{\pi}(x)} \sum_{j=1}^{m} \lambda_j [\mathbb{I}(x \in B_{2j}(x)) - \mathbb{I}(x \in B_{2j-1}(x))], 
\end{align*}
where 
$\lambda_1,\ldots,\lambda_m$ will be chosen to be uniformly distributed on $\{-1,+1\}$, and $h_1, h_2$ will be chosen to ensure
that $\varepsilon/2 \leq \pi_{\lambda}, \mu_{\lambda} \leq 1 - \varepsilon/2$.
We will set, 
\begin{align*}
\eta_\lambda(x) = \frac{\widehat{\mu}(x) - \mu_{\lambda}(x)}{1 - \pi_{\lambda}(x)}. 
\end{align*}
Now, we take a point null which we denote by $p$, to be defined by the quadruple $(\widehat{\mu},\widehat{\pi}, 0, 1)$. The functional under the null takes the value,
\begin{align*}
\int \widehat{\pi} \widehat{\mu} p_X dx = \int \widehat{\pi} \widehat{\mu}.
\end{align*}
Under the alternate, which we denote by $q_{\lambda}$ we consider the mixture defined by the quadruple $(\mu_{\lambda}, \pi_{\lambda}, \eta_{\lambda}, 1)$, and note that 
the functional takes value,
\begin{align*}
\int \pi_{\lambda} \mu_{\lambda} p_X dx \geq \int \widehat{\pi} \widehat{\mu} + m h_1 h_2,
\end{align*}
where we use the facts that the different bumps $\mathbb{I}(x \in B_{j}(x))$ do not overlap, and that $\int \mathbb{I}(x \in B_{j}(x))^2 dx = \text{vol}(B_j)$. 

By construction, $\int (\widehat{\pi}(X) - \pi_{\lambda}(X))^2 p_X(x) dx \lesssim h_1^2/\text{vol}(B_1)$, and  $\int (\widehat{\mu}(X) - \mu_{\lambda}(X))^2 p_X(x) dx \lesssim h_2^2\text{vol}(B_1)$. 
So we can choose, 
$h_1 = \sqrt{\text{vol}(B_1)}\min\{\sqrt{r_n}, \varepsilon/2\}$, $h_2 =  \sqrt{\text{vol}(B_1)} \min\{\sqrt{s_n}, \varepsilon/2\}$, to ensure that all the resulting nuisance functions are valid, and
belong to $\mathcal{G}(r_n,s_n)$.

It remains to bound the Hellinger distance for which we will use the main result of~\citet{robins2009semiparametric}. In particular, we focus on upper bounding the 
terms $a, b, d, p_j$, $j \in \{1,\ldots, m\}$, in the preamble to their Theorem 2.1, which in turn yields a bound on the squared Hellinger distance. 
We follow their notation closely. The sample space is given by $[0,1]^d \times \{0,1\} \times \{0,1\}$ which 
we partition into the sets $\mathcal{X}_j$ which are $\{0,1\} \times \{0,1\} \times B_j \cup B_{j+1}$, and so the terms $p_j$ are simply each equal to $1/m$.
Under the null, we have that,
\begin{align*}
p(X,A,Y) = \widehat{\mu}(X)^Y (1 - \widehat{\mu}(X))^{1-Y} \widehat{\pi}(X)^A (1 - \widehat{\pi}(X))^{1- A}. 
\end{align*}
On the other hand, under the alternate we have that,
\begin{align*}
q_{\lambda}(X,A,Y) &= \pi_{\lambda}(X)^A \widehat{\mu}(X)^{AY} (1 - \widehat{\mu}(X))^{A(1-Y)} (\mu_{\lambda}(X) - \pi_{\lambda}(X) \widehat{\mu}(X))^{(1-A)Y} \times \\
&~~~(1 - \pi_{\lambda}(X) - \mu_{\lambda}(X) + \pi_{\lambda}(X) \widehat{\mu}(X))^{(1-A)(1-Y)}.
\end{align*}
It is easy to verify that if we denote $q = \mathbb{E}_{\lambda}(p_{\lambda})$, then $p = q$, so the term $d$ of \citet{robins2009semiparametric} is 
0. Since we do not use a mixture under the null the term $a$ is also 0, so it only remains to bound the term $b$. 
We have that,
\begin{align*}
p - q_{\lambda} &= (1 - A)\times (-1)^Y ( \mu_{\lambda} - \widehat{\mu}) + (-1)^A \widehat{\mu}^Y (1 - \widehat{\mu})^{1-Y} (\pi_{\lambda} - \widehat{\pi}).
\end{align*}
Since, $\varepsilon/2 \leq \pi_{\lambda} \leq 1 - \varepsilon/2$ we obtain that $p > 0$. A direct calculation then gives that, the term $b$ of \citet{robins2009semiparametric} is upper bounded upto 
constants by $m (h_1^2 + h_2^2)$. Theorem 2.1 of \citet{robins2009semiparametric} 
then yields a bound on the Hellinger distance between the product measures:
\begin{align*}
H^2( p^n, \mathbb{E}_{\lambda}(q^n_{\lambda})) \lesssim m n^2 (h_1^4 + h_2^4) \lesssim \frac{n^2}{m} (r_n^2 + s_n^2).
\end{align*}
So taking $m$ sufficiently large, we obtain that $H^2( p^n, \mathbb{E}_{\lambda}(q^n_{\lambda})) \leq \alpha$, and via Lemma~\ref{lem:fuzzy} we obtain a lower bound
of order $m^2 h_1^2 h_2^2 \gtrsim r_n \times s_n$ as desired. 

\section{Proof of Theorem~\ref{thm:mainub}}
The analysis of the plugin and first-order estimators in the Gaussian sequence model follow directly from Lemma~\ref{lem:err} of the next section. We focus first on analyzing the higher-order estimator in~\eqref{eqn:hoseq}
before turning our attention to the integral of the squared density and the expected conditional covariance.

\subsection{Upper Bounds for the Quadratic Functional in the Gaussian Sequence Model}
We assume for simplicity that $y^1$ and $y^2$ are independent observations from the Gaussian sequence model. 
Observe that,
\begin{align*}
\mathbb{E} |\higherseq - \trueseq|^2 &=  |\mathbb{E} \higherseq - \trueseq|^2 + \mathbb{E} (\higherseq - \mathbb{E} \higherseq)^2 \\
&= \left[\sum_{j = T+1}^{\infty} (\widehat{\theta}_j - \theta_j^*)^2 \right]^2 +  \frac{2 \sum_{j=T+1}^\infty \widehat{\theta}_j^{2}}{n} + \frac{2 \sum_{j=1}^T \theta_j^{*2}}{n} + \frac{T}{n^2} \\
&\lesssim r_n^2 + \frac{r_n}{n} +  \frac{\|\widehat{\theta}\|_2^2}{n} + \frac{T}{n^2}, 
\end{align*}
as claimed. 
\subsection{Upper Bounds for the Integral of the Squared Density}
We first analyze the plugin estimator. 
\begin{align*}
|\plugindens - \truedens|^2 &= \left|\int \widehat{f}^2 - \int f^{*2} \right|^2 \\
&= \left|\int (\widehat{f} - f^*)(\widehat{f} + f^*)\right|^2 \\
&\leq \left|\int (\widehat{f} - f^*)^2 \right|  \left| \int (2\widehat{f} + (f^* - \widehat{f}))^2 \right| \\
&\lesssim  \left|\int (\widehat{f} - f^*)^2\right| \left( \int \widehat{f}^2 + \int  (f^* - \widehat{f}))^2\right) \\
&\leq r_n^2 + r_n \int \widehat{f}^2,
\end{align*}
as claimed. For the first-order estimator we see that,
\begin{align*}
\mathbb{E} |\firstdens - \truedens|^2 &= |\mathbb{E} \firstdens - \truedens|^2 + \mathbb{E} (\firstdens - \mathbb{E} \firstdens)^2 \\
&= \Big|\int (\widehat{f} - f^*)^2\Big|^2 + 4 \frac{\text{var}(\widehat{f}(X))}{n} \\
&\lesssim r_n^2 + \frac{\text{var}(\widehat{f}(X))}{n}. 
\end{align*}

\subsection{Upper Bounds for the Expected Conditional Covariance}
Once again, we first analyze the plugin estimator. We observe that,
\begin{align*}
\mathbb{E} |\plugincov - \truecov|^2 &= |\mathbb{E} \plugincov - \truecov|^2 + \mathbb{E} (\plugincov - \mathbb{E} \plugincov)^2 \\
&=  \left| \int \widehat{\pi} \widehat{\mu} p_X -  \int \pi^* \mu^* p_X\right|^2 + \frac{\text{var}(AY - \widehat{\pi}(X) \widehat{\mu}(X))}{n}.
\end{align*}
For the first term we observe that,
\begin{align*}
\left| \int \widehat{\pi} \widehat{\mu} p_X -  \int \pi^* \mu^* p_X\right|^2 &\lesssim \left| \int (\widehat{\mu} - \mu^*) \widehat{\pi} p_X \right|^2 + \left| \int (\widehat{\pi} - \pi^*) \mu^* p_X \right|^2 \\
&\lesssim \left(\int (\widehat{\mu} - \mu^*)^2 p_X \right) \left( \int \widehat{\pi}^2 p_X \right) +  \left(\int (\widehat{\pi} - \pi^*)^2 p_X \right)  \left( \int \mu^{*2} p_X \right) \\
&\lesssim r_n  \left( \int \widehat{\pi}^2 p_X \right) + s_n \left( \int (\widehat{\mu}^{2} + (\widehat{\mu} - \mu^*)^2)  p_X \right) \\
&\lesssim r_n \times s_n + r_n \left( \int \widehat{\pi}^2 p_X \right) + s_n \left( \int \widehat{\mu}^2 p_X \right),
\end{align*}
as desired. For the first-order estimator we have,
\begin{align*}
\mathbb{E} |\firstcov - \truecov|^2 &= |\mathbb{E} \firstcov - \truecov|^2 + \mathbb{E} (\firstcov - \mathbb{E} \plugincov)^2 \\
&= \left| \int (\pi^* - \widehat{\pi})(\mu^* - \widehat{\mu}) p_X \right|^2 +  \frac{\text{var}((A - \widehat{\pi}(X))(Y - \widehat{\mu}(X)))}{n} \\
&\leq \left(\int (\pi^* - \widehat{\pi})^2 p_X \right) \left( \int (\mu^* - \widehat{\mu})^2 p_X \right) +  \frac{\text{var}((A - \widehat{\pi}(X))(Y - \widehat{\mu}(X)))}{n} \\
&\leq r_n \times s_n +  \frac{\text{var}((A - \widehat{\pi}(X))(Y - \widehat{\mu}(X)))}{n},
\end{align*}
as claimed. 

\section{An Adaptive Estimate of the Quadratic Functional}
\label{app:super}

Throughout our paper we focused our discussion on the case when 
the error of the pilot estimate was larger than the variance of the first-order estimator. This situation is typical. In this section, we briefly investigate
the setting where the pilot may be super-accurate. 

\subsection{Upper Bounds} 
For any $\delta > 0$, consider the following adaptive estimate:
\begin{align*}
\adaseq = \begin{cases}
\pluginseq~\text{if}~|\pluginseq - \firstseq| \leq  4 \|\widehat{\theta}\|_2 \sqrt{\frac{4\log(2/\delta)}{n}}\\
\firstseq~\text{otherwise}.
\end{cases}
\end{align*}
The estimate is similar to a Lepski-style adaptive estimator which chooses between the plugin and first-order estimate. The following result then holds:
\begin{theorem}
\label{thm:adaub}
For any $\delta > 0$, the risk of $\adaseq$ is upper bounded as:
\begin{align*}
\mathbb{E}[(\adaseq - Q(\theta^*))^2] \lesssim r_n^2 + \min \left\{r_n \|\widehat{\theta}\|_2^2 +  \frac{ \delta \|\widehat{\theta}\|_2^2}{n},  \frac{\|\widehat{\theta}\|^2_2 \log(1/\delta)}{n} \right\}. 
\end{align*}
\end{theorem}

\noindent {\bf Remarks: } 
\begin{enumerate}
\item Suppose that we ignored the terms which depend on $\delta$, then the estimator $\adaseq$ is (fully) adaptive, i.e. it achieves the same performance as an oracle which always picked the estimate 
with lower risk. 
\item More generally, there is a small price to pay to adapt between these estimates. This price is closely related to the standard Hodges superefficiency phenomenon in estimating a Gaussian mean.
We develop this further in Appendix~\ref{app:adalb}.
\end{enumerate}

\begin{proof}
Let us define $\|\widehat{\theta} - \theta^*\|^2_2 := R$.
We first note the following error bounds on our estimates:
\begin{lemma}
\label{lem:err}
We have the following bounds:
\begin{align*}
|\pluginseq - Q(\theta^*)|^2 &\lesssim R^2 + R \|\widehat{\theta}\|_2^2 \\
\mathbb{E}|\firstseq - Q(\theta^*)|^2 &\lesssim R^2 + \frac{ \|\widehat{\theta}\|_2^2}{n} \\
\sqrt{\mathbb{E}|\firstseq - Q(\theta^*)|^4} &\lesssim R^2 + \frac{ \|\widehat{\theta}\|_2^2}{n} .
\end{align*}
\end{lemma}
Taking this result as given we can complete the proof. 
Recall, that we observe $y = \theta^* + \epsilon$. We define the event $E$ to be the event on which, 
\begin{align*}
|\inprod{\epsilon}{\widehat{\theta}}| \leq \|\widehat{\theta}\|_2 \sqrt{\frac{4 \log (2/\delta)}{n}},
\end{align*}
which happens with probability at least $1 - \delta^2$, by applying a standard Gaussian tail bound.

Now observe that when $R  \leq \frac{4 \log (2/\delta)}{n}$, and on the event $E$:
\begin{align*}
|\pluginseq - \firstseq| &= |2 (\|\widehat{\theta}\|_2^2 - \inprod{y}{\widehat{\theta}})| \leq 2 \|\widehat{\theta}\|_2 \sqrt{\frac{4 \log (2/\delta)}{n}} + 2 |\inprod{\widehat{\theta}}{\theta^* - \widehat{\theta}}| \\
&\leq 2 \|\widehat{\theta}\|_2 \sqrt{\frac{4 \log (2/\delta)}{n}} + 2 \sqrt{R} \|\widehat{\theta}\|_2 \leq 4 \|\widehat{\theta}\|_2 \sqrt{\frac{4 \log (2/\delta)}{n}},
\end{align*}
so our selection rule picks the estimate $\pluginseq$.

Now, we are in a position to analyze our selection rule. 
Let us denote an index $j$ which takes value $1$ when $\adaseq = \pluginseq$ and $2$ otherwise. 
We consider two cases, when $R  \leq  \sqrt{\frac{4 \log (2/\delta)}{n}}$ and when $R >  \sqrt{\frac{4 \log (2/\delta)}{n}}$.
In the first case,
\begin{align*}
\mathbb{E}[(\adaseq - Q(\theta^*))^2] &= \mathbb{E}[(\adaseq - Q(\theta^*))^2 \mathbb{I}[j = 1]] + \mathbb{E}[(\adaseq - Q(\theta^*))^2 \mathbb{I}[j = 2]] \\
&\lesssim R^2 + R \|\widehat{\theta}\|_2^2 + \delta \sqrt{ \mathbb{E}[(\adaseq - Q(\theta^*))^4} \\
&\lesssim R^2 + R \|\widehat{\theta}\|_2^2 + \delta \left[ R^2 + \frac{ \|\widehat{\theta}\|_2^2}{n}\right] \\
&\lesssim R^2 + R \|\widehat{\theta}\|_2^2 +  \frac{ \delta \|\widehat{\theta}\|_2^2}{n}.
\end{align*}
In the second case, 
\begin{align*}
\mathbb{E}[(\adaseq - Q(\theta^*))^2] &= \mathbb{E}[(\adaseq - Q(\theta^*))^2 \mathbb{I}[j = 1]] + \mathbb{E}[(\adaseq - Q(\theta^*))^2 \mathbb{I}[j = 2]] \\
&= \mathbb{E}[(\pluginseq - Q(\theta^*))^2 \mathbb{I}[j = 1]] + \mathbb{E}[(\firstseq - Q(\theta^*))^2 \mathbb{I}[j = 2]] \\
&\lesssim \mathbb{E}[(\pluginseq - \firstseq)^2 \mathbb{I}[j = 1]] +  \mathbb{E}[(\firstseq - Q(\theta^*))^2] \\
&\lesssim  R^2 + \frac{\|\widehat{\theta}\|^2_2 \log(1/\delta)}{n}.
\end{align*}
Putting these together we obtain the desired theorem.
\end{proof}

\noindent {\bf Proof of Lemma~\ref{lem:err}: } We prove each of the three claims in turn. 
Observe that,
\begin{align*}
|\pluginseq - Q(\theta^*)|^2 &= |\|\widehat{\theta}\|_2^2 - \|\theta^*\|_2^2|^2 \\
&= |(\widehat{\theta} - \theta^*)^T (\widehat{\theta} + \theta^*)|^2 \\
&\lesssim R (\|\widehat{\theta}\|_2^2 + R), 
\end{align*}
as desired. Now, we see that,
\begin{align*}
\mathbb{E}|\firstseq - Q(\theta^*)|^2 &= (\mathbb{E}(\firstseq) - Q(\theta^*))^2 +\mathbb{E}( \mathbb{E}(\firstseq) - \firstseq)^2 \\
&= \|\widehat{\theta} - \theta^*\|_2^4 + \frac{4 \|\widehat{\theta}\|_2^2}{n} \\
&\lesssim R^2 +  \frac{ \|\widehat{\theta}\|_2^2}{n},
\end{align*}
where the bounds on the bias and variance follow from a direct calculation.
Finally,
\begin{align*}
\sqrt{\mathbb{E}|\firstseq - Q(\theta^*)|^4} &\lesssim \sqrt{ (\mathbb{E}(\firstseq) - Q(\theta^*))^4 + \mathbb{E}( \mathbb{E}(\firstseq) - \firstseq)^4} \\
&\lesssim R^2 + \frac{ \|\widehat{\theta}\|_2^2}{n},
\end{align*}
where to bound the second term, we simply use the fact that the fourth moment of a mean zero Gaussian random variable is $3 \sigma^4$. 

\subsection{An Adaptivity Lower Bound}
\label{app:adalb}
In this section, we prove the following complementary lower bound which shows a sense in which our adaptive estimator is unimprovable. 
We denote the risk of our adaptive estimate $\adaseq$ as:
\begin{align*}
f_{\delta}(r) = r^2 + \min \left\{ r \|\widehat{\theta}\|_2^2 + \frac{\delta \|\widehat{\theta}\|_2^2}{n}, \frac{\log(1/\delta) \|\widehat{\theta}\|_2^2}{n}\right\},
\end{align*}
where $\delta > 0$ is chosen to be sufficiently small. 
\begin{lemma}
\label{lem:adalb}
Suppose that we have an estimate $\widehat{Q}$ such that for some $r_1 \geq 0$, and for a sufficiently small $\varepsilon > 0$, 
\begin{align*}
\sup_{\theta^* \in \Theta(r_1)} \mathbb{E}(\widehat{Q} - Q(\theta^*))^2 \lesssim \varepsilon f_{\delta}(r_1)
\end{align*}
then there exists an $r_2$, such that,
\begin{align*}
\sup_{\theta^* \in \Theta(r_2)} \mathbb{E}(\widehat{Q} - Q(\theta^*))^2 \gtrsim \frac{\log(1/(\varepsilon \delta))}{\log(1/\delta)} f_{\delta}(r_2),
\end{align*}
so long as $\|\widehat{\theta}\|^2_2 \geq \log^2(1/(\epsilon \delta))/(n \log(1/\delta))$. 
\end{lemma}
\noindent This result captures the fact that there is a strong sense in which Theorem~\ref{thm:adaub} achieves the limits of adaptation.
In particular, if an estimator has a (very) small risk relative to $\adaseq$ for some value of the unknown radius $r_{1}$, then this must come
at the expense of a worse performance at a different value $r_{2}$. 
Our lower bound follows the well-trodden route of using the
constrained risk inequality in Lemma~\ref{lem:cri} to argue that achieving a (very) small risk at a point in the parameter space comes at the expense of a larger risk in the neighborhood of that point.

We note the (mild) restriction that $\|\widehat{\theta}\|^2_2 \geq \log^2(1/(\epsilon \delta))/(n \log(1/\delta))$. When this assumption does not hold, the adaptation picture is different. This can be seen by observing that when $\|\widehat{\theta}\|_2 = 0$, 
the estimate $\widehat{Q} = 0$ is adaptively optimal for any $r$, and there is no price to pay for adaptation.
\begin{proof}

We consider two cases, when $r_1 \leq \log(1/\delta)/n$ and when $r_1 \geq \log(1/\delta)/n$. 

\noindent {\bf Case 1: } When $r_1 \leq \log(1/\delta)/n$ we can write,
\begin{align*}
f_{\delta}(r_1) \leq r_1^2 + r_1 \|\widehat{\theta}\|_2^2 + \frac{\delta \|\widehat{\theta}\|_2^2}{n}, 
\end{align*}
since $\delta > 0$ is sufficiently small to ensure that $\delta \leq \log(1/\delta)$. Now, this implies that $\widehat{Q}$ satisfies,
\begin{align}
\label{eqn:temp}
\sup_{\theta^* \in \Theta(r_1)} \mathbb{E}(\widehat{Q} - Q(\theta^*))^2 \lesssim \varepsilon \left(r_1^2 + r_1 \|\widehat{\theta}\|_2^2 + \frac{\delta \|\widehat{\theta}\|_2^2}{n}\right).
\end{align}
The lower bounds in Appendix~\ref{app:seqlb} already show that, 
\begin{align*}
\sup_{\theta^* \in \Theta(r_1)} \mathbb{E}(\widehat{Q} - Q(\theta^*))^2 \gtrsim r_1^2 + r_1 \|\widehat{\theta}\|_2^2,
\end{align*}
so the claimed improvement is only possible if the final term in~\eqref{eqn:temp} dominates. In this case we have that,
\begin{align*}
\sup_{\theta^* \in \Theta(r_1)} \mathbb{E}(\widehat{Q} - Q(\theta^*))^2 \lesssim  \frac{\varepsilon \delta \|\widehat{\theta}\|_2^2}{n}.
\end{align*}
We now apply Lemma~\ref{lem:cri} with the choices, $\theta := \widehat{\theta}$, $\beta := \sqrt{\varepsilon \delta} \|\widehat{\theta}\|_2/\sqrt{n}$, and 
$\alpha := \sqrt{\log(1/(\varepsilon \delta))/n}$, to conclude that, 
\begin{align*}
\sup_{\theta^* \in \Theta(\alpha^2)} \mathbb{E}(\widehat{Q} - Q(\theta^*))^2 \gtrsim \alpha^4 + \frac{\log(1/(\varepsilon \delta)) \|\widehat{\theta}\|_2^2}{n} \gtrsim \frac{\log(1/(\varepsilon \delta)) \|\widehat{\theta}\|_2^2}{n}.
\end{align*}
On the other hand, we have that the risk of $\adaseq$ is upper bounded as:
\begin{align*}
f_{\delta}(\alpha^2) \leq \alpha^4 + \frac{\log(1/\delta) \|\widehat{\theta}\|_2^2}{n} \lesssim  \frac{\log(1/\delta) \|\widehat{\theta}\|_2^2}{n},
\end{align*}
using our assumed lower bound on $\|\widehat{\theta}\|_2$. This in turn establishes the non-adaptivity claim. 

\vspace{.3cm}

\noindent {\bf Case 2: } We now consider the case when $r_1 \geq  \log(1/\delta)/n$. 
The lower bounds in Appendix~\ref{app:seqlb} already show that 
\begin{align*}
\sup_{\theta^* \in \Theta(r_1)} \mathbb{E}(\widehat{Q} - Q(\theta^*))^2 \gtrsim r_1^2 + \frac{\|\widehat{\theta}\|_2^2}{n},
\end{align*}
so once again the claimed improvement is only possible if the second term dominates and we have,
\begin{align*}
\sup_{\theta^* \in \Theta(r_1)} \mathbb{E}(\widehat{Q} - Q(\theta^*))^2 \lesssim \frac{\varepsilon \log(1/\delta) \|\widehat{\theta}\|_2^2}{n}.
\end{align*}
Now, once again we apply Lemma~\ref{lem:cri} with the choices, $\theta := \widehat{\theta}$, $\beta := \sqrt{\varepsilon \log(1/\delta)} \|\widehat{\theta}\|_2/\sqrt{n}$,
and $\alpha := \sqrt{\log(1/(\varepsilon \log(1/\delta))/n}$. To obtain that,
\begin{align*}
\sup_{\theta^* \in \Theta(0)} \mathbb{E}(\widehat{Q} - Q(\theta^*))^2 \gtrsim \frac{\log(1/(\varepsilon \log(1/\delta))) \|\widehat{\theta}\|_2^2}{n}.
\end{align*}
On the other hand, the estimator $\adaseq$ achieves the guarantee, 
\begin{align*}
\sup_{\theta^* \in \Theta(0)} \mathbb{E}(\adaseq - Q(\theta^*))^2 \lesssim \frac{\delta \|\widehat{\theta}\|_2^2}{n}.
\end{align*}
So we obtain that, taking $r_2 = 0$, 
\begin{align*}
\sup_{\theta^* \in \Theta(r_2)} \mathbb{E}(\widehat{Q} - Q(\theta^*))^2 \gtrsim \frac{\log(1/(\varepsilon \log(1/\delta)))}{\delta} f_{\delta}(r_2) \gtrsim \frac{\log(1/(\varepsilon \delta))}{\log(1/\delta)} f_{\delta}(r_2).
\end{align*}
These two facts taken together yield the non-adaptivity claim of the theorem.

\end{proof}

\end{document}